\documentclass[a4paper,10pt,oneside,onecolumn]{article}

\usepackage{geometry}
\usepackage{graphicx}%
\usepackage{labelfig}%
\usepackage{diagbox}%

\usepackage{amsmath,amsfonts,amssymb,amsthm,upref,bm}%

\usepackage{mathrsfs}%
\usepackage{fancyvrb}%
\usepackage{upgreek}%

\makeatletter
\renewcommand\section{\@startsection {section}{1}{\z@}%
           {18\p@ \@plus 6\p@ \@minus 3\p@}%
           {9\p@ \@plus 6\p@ \@minus 3\p@}%
           {\normalsize\bfseries\boldmath}}

\DeclareMathAlphabet{\varmathbb}{U}{pxsyb}{m}{n}

\newcommand{\MF}[1]{\mathop{#1\vrule height0.45ex width0pt}\nolimits}

\newcommand{\pd}[3][]{\mathchoice{\raise-0.5pt\hbox{$\partial$}%
\vphantom{\partial}_{\mkern-1.5mu#2}^{\mkern0.4mu#1}\mkern0.3mu}%
{\raise-0.5pt\hbox{$\partial$}%
\vphantom{\partial}_{\mkern-1.5mu#2}^{\mkern0.4mu#1}\mkern0.3mu}%
{\raise-0.5pt\hbox{$\scriptstyle\partial$}%
\vphantom{\partial}_{\mkern-1.7mu#2}^{\mkern0.1mu#1}\mkern0.1mu}%
{\raise-0.5pt\hbox{$\scriptscriptstyle\partial$}%
\vphantom{\partial}_{\mkern-1.7mu#2}^{\mkern0.1mu#1}\mkern0.1mu}#3}

\newcommand{\ii}{\kern0.05em\mathrm{i}\kern0.05em} \newcommand{\E}[1]{\textrm{e}^{#1}}
\newcommand{\D}{\mathrm{d}\kern0.2pt} \renewcommand{\vec}[1]{\bm{#1}}

\renewcommand{\Re}{\MF{\mathrm{Re}}} \renewcommand{\Im}{\MF{\mathrm{Im}}}

\newcommand{\RR}{\varmathbb{R}} 
\newcommand{\transp}{\textsf{T}}

\newcommand{\oscf}{\varpi}
\newcommand{\accuracy}{\varepsilon}

\allowdisplaybreaks[3]

\begin{document}

\begin{center}

{\Large On computation of oscillating integrals of ship-wave theory}

\vskip6mm

{\large Oleg\ V.\ Motygin}

\vskip3mm

{\itshape Institute of Problems in Mechanical Engineering, Russian Academy\\ of
Sciences, V.O., Bol'shoy pr., 61, 199178 St.\,Petersburg, Russia}\\
e-mail: \texttt{o.v.motygin@gmail.com}

\end{center}

\vskip3.5mm

\noindent{\bfseries Abstract}

\vskip2mm

\noindent Green's function of the problem describing steady forward motion of
bodies in an open ocean in the framework of the linear surface wave theory (the
function is often referred to as Kelvin's wave source potential) is considered.
Methods for numerical evaluation of the so-called `single integral' (or, in
other words, `wavelike') term, dominating in the representation of Green's
function in the far field, are developed.  The difficulty in the numerical
evaluation is due to integration over infinite interval of the function
containing two differently oscillating factors and the presence of stationary
points. This work suggests two methods to approximate the integral and its
derivatives. First of the methods is based on the idea suggested by D.\,Levin in
1982~---~evaluation of the integral is converted to finding a particular slowly
oscillating solution of an ordinary differential equation. To overcome
well-known numerical instability of Levin's collocation method, an alternative
type of collocation is used; it is based on a barycentric Lagrange interpolation
with a clustered set of nodes. The second method for evaluation of the wavelike
term involves application of the steepest descent method and Clenshaw--Curtis
quadrature. The methods are numerically tested and compared for a wide variety
of parameters.

\vskip4mm

\noindent\textit{Keywords: } 

{\raggedright

\vskip1mm

Ship waves, Kelvin wave source, Green's function, Oscillatory integral, Levin method,
Clenshaw--Curtis quadrature, Barycentric Lagrange interpolation, Chebyshev points


\vskip1mm


}



\section{\normalsize Introduction}
\label{intro0}

The steady motion of a ship in calm water is a classic problem of hydrodynamics,
significant for engineering and ship design because of its applicability to
prediction of wave pattern and wave resistance. As early as in 1887, Lord Kelvin
\cite{Kelvin1887} first provided a mathematical description of the wave pattern
created by an object moving forward with a constant speed; the familiar V-shaped
pattern behind a ship (still being in study, see e.g.\
\cite{ReedMilgram2002,RabaudMoisy2013,DarmonBenzaquenRaphael2014} and references
therein) now bears his name. Starting from Lord Kelvin's work, substantial
efforts have been spent to investigation of the ship wave problem and
determining the wave resistance. At this, an extensive literature has been
produced; comprehensive reviews can be found, e.g.\ in
\cite{WehausenLaitone,Wehausen73,Miloh91,Eggers,LarssonBaba96,LWW}. In
particular, much attention has been paid to potential models and linearized
statements similar to that considered in the present work.

An important role in the linear theory of ship waves and wave resistance is
played by the Kelvin wave source potential, which may also be identified as the
Green's function of the Neumann--Kelvin problem. This potential describes a
source moving with a constant horizontal velocity in an ideal incompressible
fluid having a free surface; fluid's motion is irrotational and steady-state in
a coordinate system attached to the source. The Kelvin wave source potential is
fundamental in the theory and applications~---~for instance, a solution of the
problem of the flow past a moving ship is often sought as a distribution of the
sources on ship's hull. Considerable efforts in the analysis of the Green's
function has resulted in a vast literature (see \cite{LWW} and references
therein) with special interest to the far field behaviour of the source (see
e.g.\ \cite{Eggers,MV}) and to development of efficient and accurate algorithms
for computing Kelvin's potential.

A number of representations of the Green's function have been derived (see
\cite{Noblesse}). Generally, the Green's function is divided into three terms:
the singular (Rankine) one, whose computation is elementary; the wavelike
(referred to as `single integral') term~---~oscillating and dominating in the
far field; the near-field nonoscillatory term (for historical reasons often
referred to as `double integral term'). Detailed mathematical analysis of the
behavior of the near-field and the far-field terms can be found in
\cite{Euvrard83,MV,LWW}. It is known that the `double integral term' can be
expressed as a single integral, in particular, a number of such representations
in terms of the exponential integral function can be found in
\cite{ShenFarell1977,Noblesse}. This allows realization of fast and effective
schemes for numerical evaluation of the `double integral term'. A practical
numerical implementation can be found in \cite{Newman1}; further development for
the efficacy and accuracy in computation of the near-field term is described in
\cite{PonizyEtAl1998}.

At the same time, the question of evaluation of the wavelike term is only partly
solved. The main difficulty in computation of the term is due to the presence of
two oscillating factors of different nature in the integrand, which can have
stationary points, and the infinite interval of integration. Besides, there is a
difficult limiting case: it was shown in \cite{Ursell60} (more details are given
in \cite{Euvrard83,Ursell88}) that the track of the source moving in the free
surface is a line of essential singularities. Near the line the wave elevation
oscillates with indefinitely increasing amplitude and indefinitely decreasing
wavelength. (It is interesting that introducing of surface tension in the
formulation of the problem is reported in \cite{Chen2002} to eliminate the
singularity of the Green's function.)

Most of the existing methods for evaluation of the wavelike term are based or
involve two very useful expansions given in \cite{Bessho64}: convergent series
and an asymptotic one (completed in \cite{Ursell60}). However, the summation of
such series is not a trivial computational task, it leads to poor accuracy when
the source and the field points are close to the free-surface, the method is
vulnerable to severe numerical problems because of the presence of very large in
magnitude alternating terms and the criteria to choose one of the expansions is
difficult to determine. Details of numerical algorithms using the expansions,
discussion of their applicability and heuristic criteria of choice between the
algorithms can be found in \cite{BaarPrice,BaarPrice2,Marr,PonizyEtAl1998}.
Another approach in \cite{WangRogers1990} to calculation of the wavelike term is
also based on the Green's function representation \cite{Bessho64} but the work
is in terms of integrals as opposed to the use of the series. Computation of the
wavelike term near the track was discussed in \cite{Newman3} where expansions of
the function and numerical algorithms are given. 

In the present work our purpose is to elaborate alternative, accurate and fast
computation methods to approximate the wavelike term. Our consideration is based
on recently developed techniques for evaluation of integrals of oscillating
functions. The quadrature of oscillating integrals is a very important
computational problem (widely considered as a difficult one) appearing in a wide
range of applications. The field has enjoyed a recent upsurge of interest and
substantive developments with a number of old methods being enhanced and new
ones being suggested. Amidst them we mention the numerical steepest descent (see
\cite{HuybrechsVandewalle2007,DeanoHuybrechs2009}), Filon-type (see
\cite{Filon1928,Flinn1960,IserlesNorsett2005,Olver2007,Xiang2007b}), Levin-type
methods (discussed below), Filon--Clenshaw--Curtis quadrature (see
\cite{EvansWebster1999,DominguezGrahamKim2013}), `double exponential formula'
for Fourier-type integrals (see \cite{Ooura2005} and references therein).
Reviews of the methods, their comparison and analysis, as well as bibliography
notes can be found in
\cite{EvansWebster1999,INO06,Olver2008,HuybrechsOlver2009}.

The first of two schemes developed in the present paper is based on the ideas
of Levin, suggested in \cite{Levin82,Levin97}. The Levin method converts the
calculation of an oscillatory integral into solving an ordinary differential
equation (ODE) whose special, slowly oscillating solution is sought. Typically,
the Chebyshev expansion is used as a representation of the solution and a
collocation reduces the ODE to a system of linear equations. The Levin method
has attracted much attention due to its ability to handle integrals with
complicated phase functions, but the collocation has been found to be
numerically unstable, leading to an ill-conditioned linear system for large
number of nodes. Analysis and efforts to improve the stability of the Levin
collocation method have been made, a number of Levin-type methods have been
developed (see
\cite{EvansWebster1999,EvansChung2003,INO06,Xiang2007b,LiWangWang2008,Olver2010,Li2010,Li2011}
and references therein). Some progress has been achieved using careful analysis
of the linear system and advanced methods of linear algebra such as TSVD
(truncated singular value decomposition) and GMRES (generalized minimal
residual) method. Our scheme is based on a different type of
collocation~---~instead of the Chebyshev expansion we use a special form of
Lagrange polynomial interpolation which was recently discovered as an effective
tool for numerical analysis.

It is to note that polynomial interpolation formulae are widely used in
theoretical studies, but not in numerical practice~---~in view of their instability.
Finding the polynomials involves solution of a Vandermonde linear system of
equations, which is exponentially unstable. Runge's phenomenon is also well
known: for equispaced interpolation points small perturbations of the initial
data may result in huge changes of the interpolant. However, it has been found
(see \cite{Salzer1972,Henrici1979,BerrutTrefethen2004,TrefethenBook} and
references therein) that these problems are avoided when using the Lagrange
interpolation in one of the so-called barycentric forms and with the nodes
clustered near the ends of the interval of approximation, e.g.\ at Chebyshev
points. A review of the existing results and explanation of the attractive
features of the barycentric Lagrange interpolation with a clustered set of nodes
can be found in \cite{BerrutTrefethen2004,TrefethenBook}; rigorous confirmation
of stability is given in \cite{Higham}.

Thus, in our first scheme, following \cite{Levin82,Levin97}, we reduce
evaluation of the integral in question to solution of an ODE on a finite
interval. We prove that the ODE has one bounded solution and define the value of
the function at the right end of the interval. Solution's value at the left end
of the interval, up to a known factor, coincides with the sought value of the
integral. To find the solution of the ODE numerically, we seek it in the form of
a barycentric Lagrange interpolating polynomial (alternatively, the
representation can include a term arising from an asymptotic analysis, intended
to absorb `bad' features of the solution such as sharp peaks). Then we apply
collocation of the equation on a set of Chebyshev points. Our numerical
experiments show that the constructed numerical algorithm is stable, typically
converging steadily to the level of rounding errors with increase of the number
of nodes, even in the presence of stationary points (which is considered as
complication for Levin and Levin-type methods).

The second scheme, developed in the present paper, relies upon Clenshaw--Curtis
quadrature \cite{ClCu} and the steepest descent method. First, we transform the
integral contour along the steepest descent path in the complex plane and change
variable to have a finite interval of integration. Further computation is based
on the Clenshaw--Curtis quadrature, the method for numerical integration derived
by an expansion of the integrand in terms of Chebyshev polynomials and
approximation of the expansion coefficients using discrete Fourier
transformation. The advantages of the quadrature rule are its fast convergence
for smooth integrands at increasing number of nodes $n$ (see
\cite{OHaraSmith1968,Trefethen,TrefethenBook}), naturally nested quadrature
rules, effective computation of weights based on Fast Fourier Transform (see
\cite{Gentleman1972I,Waldvogel2006}). So, the Clenshaw--Curtis quadrature shares
the insensitivity of the FFT to rounding errors (see e.g.\ \cite{KanekoLiu1970})
and the cost of implementation at $\MF{O}( n \log n)$ operations. These
properties of the quadrature, the simplicity and smoothness of the considered
integrand suggest that in our particular case a `brute-force' application of the
quadrature can be fairly effective, which is confirmed by our numerical
experiments. In the suggested scheme we build a sequence of approximations by
increasing the number of quadrature's nodes (so that the number of intervals
between nodes is doubled at each step) until some last values of the sequence
become closer to each other than the given tolerance. Such simple estimate of the
approximation error is known to be very realistic for the Clenshaw--Curtis
integration (see \cite{Gentleman1972I}), which is also confirmed by our
computations.

Using these two numerical algorithms we carry out evaluation of the wavelike
term of the Green's function for a wide variety of parameters, comparing the
accuracy of the approximations and their efficacy. These numerical experiments
allow us to conclude that both developed algorithms are consistent, reliable and
compatible in speed while the algorithm based on Levin ODE and barycentric
Lagrange interpolation is somewhat faster when the order of interpolating
polynomial, needed for achieving the given accuracy, is small (e.g.\ when the
source or the field point is located sufficiently deeply). At the same time, the
algorithm based on Clenshaw--Curtis quadrature works surprisingly well in the
most difficult for the numerical integration zone, near the track of the source
moving closely to the free surface.

We also consider application of the suggested methods to computation of
derivatives of the wavelike term of the Green's function~---~the ability to
compute the derivatives is fundamental in realization of many numerical methods
(in particular, those based on integral equations of the potential theory). We
show that generalization of both suggested numerical schemes is rather
straightforward.

Now we give a brief outline of the paper. In \S\,\ref{sect:1} we present the
mathematical problem describing forward motion of a source and introduce a
representation of the Green's function~---~the following sections are devoted to
evaluation of its wavelike component. In \S\,\ref{sect:3} we reduce evaluation
of the integral to finding a special solution of an ODE and study properties of
the solution. In \S\,\ref{sect:5} a numerical scheme based on the barycentric
Lagrange interpolation is constructed for finding the solution of the
differential equation. An alternative scheme for evaluation of the wavelike
term, based on the steepest descent method and the Clenshaw--Curtis quadrature,
is developed in \S\,\ref{sect:7}. In \S\,\ref{sect:8} numerical experiments are
carried out to test and compare the algorithms. Application of the suggested
methods to computation of derivatives of the wavelike term is discussed in
\S\,\ref{sect:der}. The paper concludes with a summary of the principal
findings.

\section{Statement of the problem and a representation of the Green function}
\label{sect:1}

Consider the mathematical problem for a velocity potential that describes the
forward motion of a source with a constant speed $U$ through an unbounded heavy
fluid having a free surface $F$. The fluid is assumed to be ideal incompressible
and having infinite depth; fluid's motion is irrotational and steady-state in a
Cartesian coordinate system $(x,y,z)$ attached to the source. The $x$-axis is
parallel to the direction of motion, $z$ is another horizontal axis, and
$y$-axis is vertical and directed upwards. We will also use the notation
$\vec{x}=(x,y,z)$. The source is located at $\vec{x}_0=(x_0,y_0,z_0)$ so that
the fluid fills in $W=\RR_-^3\setminus\{\vec{x}_0\}$, where
$\RR_-^3=\{(x,z)\in\RR^2,\,y<0\}$, and $F=\partial \RR_-^3$.

We consider the mathematical problem of the linear surface wave theory (often
referred to as the Neu\-mann--Kelvin problem~---~see book \cite{LWW} and
references therein). The motion of the source is described by a velocity
potential $\MF{G}(\vec{x};\vec{x}_0)=\MF{G}(x,y,z;x_0,y_0,z_0)$
\begin{gather}
 -\nabla^2G=\MF{\updelta}(x-x_0)\MF{\updelta}(y-y_0)\MF{\updelta}(z-z_0)\quad\mbox{when}
  \quad y\leq0,\ \ y_0<0,\label{eq:1}\\
 \pd[2]{x}G+\nu\pd{y}G=0\quad\mbox{when}
  \quad y=0,\ \ y_0<0,\label{eq:2}
\end{gather}
where $\updelta$ is Dirac's delta-function,
$\nabla=(\pd{x}{},\pd{y}{},\pd{z}{})$, $\nu=g/U^2$ is the wave number, $g$
denotes acceleration due to gravity.

Since the fluid domain is unbounded, the basic relations must be complemented
with conditions prescribing behaviour of the potential at infinity (see
\cite{Gamst79,Euvrard83,LWW,M2012}). Following \cite{M2012}, we demand
\begin{equation}
 \sup_{\RR_-^3\setminus\{\vec{x}:|\vec{x}-\vec{x}_0|<r\}}|\nabla G|<\infty\quad\mbox{for some}\quad r>0,\label{eq:3}
\end{equation}
and
\begin{equation}
\iiint_{\{y<0,\,x> x_*\}}|\nabla G|^2\,\D{}x\,\D{}y\,\D{}z<\infty\quad\mbox{for
some}\quad x_*>x_0.
\label{eq:4}
\end{equation}
The latter condition is one of a few variants suggested in \cite{M2012}. The
statement \eqref{eq:1}--\eqref{eq:4} is proved in \cite{M2012} to have a unique
(up to an additive constant) solution.

Without loss of generality we can assume that $x_0=0$ and $z_0=0$. Further on we will
assign the values and use the notation:
\begin{equation*}
  \MF{G}(x,y,z;y_0)=\MF{G}(x,y,z;0,y_0,0).
\end{equation*}
Besides, further we will use dimensionless coordinates by fixing $\nu=1$.

A number of representations of the Green's function are known (see e.g.\
\cite{WehausenLaitone,LWW} and references therein). Typically, they involve one
double integral and one single integral. The single integral represents the
far-field oscillatory part where one may recognize the well-known ship wake
pattern, whilst the double integral mostly contributes to the local flow near
the source. More suitable for numerical applications representations with the
double integrals being transformed to single ones, have been derived (see
\cite{ShenFarell1977,Noblesse} and references therein). Further we use the
formula suggested in \cite{Noblesse}:
\[
  \MF{G}(x,y,z;y_0)=-\frac{1}{4\pi R}+\frac{1}{4\pi R'}+
 \MF{I_0}(x,y+y_0,z)+\MF{I_\infty}(x,y+y_0,z).
\]
Here $R=\sqrt{x^2+(y-y_0)^2+z^2}$, $R'=\sqrt{x^2+(y+y_0)^2+z^2}$, the `near-field'
term is written as follows:
\begin{equation*}
  \MF{I_0}(x,y,z)=\frac{1}{2\pi^2}
 \Im\int_{-1}^{1}\E{\MF{\xi}(x,y,z,t)}\MF{E_1}(\MF{\xi}(x,y,z,t))\,\D{}t,
\end{equation*}
where
$\MF{\xi}(x,y,z,t)=y\bigl(1-t^2\bigr)+\bigl(zt+\ii|x|\bigr)\bigl(1-t^2\bigr)^{1/2}$
and $E_1$ is the integral exponent function (see e.g.\ \cite[\S\,5.1]{AS}). The
integral representing the oscillating wavelike part of the Green function has the
following form:
\begin{equation}
  \MF{I_\infty}(x,y,z)=\pi^{-1}\MF{H}(-x)\Im\bigl\{\MF{I}(x,y,z)+\MF{I}(x,y,-z)\bigr\}.
\label{eq:defI_inf}
\end{equation}
Here $H$ is Heaviside's function,
\begin{gather}
  \MF{I}(x,y,z)=\int_0^{\infty}\E{\MF{\oscf}(t,x,y,z)}\,\D{}t,\ \ \textrm{where}\ \ y<0,\ x\leq0,\ z\in\RR,
\label{eq:defI}
\\
\MF{\oscf}(t,x,y,z)=y\left(1+t^2\right)+\ii(x+zt)\sqrt{1+t^2}.\notag
\end{gather}

Numerical evaluation of the function $\MF{I_0}(x,y,z)$ is discussed in
\cite{Noblesse1977,Noblesse,Newman1}. In the present work we are concerned with
approximation of $\MF{I}(x,y,z)$. Except the case $x=z=0$ when
$I=\frac{\sqrt\pi}{2}\frac{\E{y}}{\sqrt{-y}}$, the problem is rather difficult
because of the presence of two differently oscillating factors and the infinite
interval of integration. For $z>0$ and $x\leq-2\sqrt{2}z$ there are also
critical points
$t_\textrm{crit}^\pm=-\frac{x}{4z}\pm\sqrt{\left(\frac{x}{4z}\right)^2-\frac12}$
defined by the equation $2zt^2+xt+z=0$ ($\Im\{\oscf'_t\}=0$), merging to one
point when $x=-2\sqrt{2}z$.

\section[Evaluation of the integral $\MF{I}(x,y,z)$ based on an ordinary
differential equation]{Evaluation of the integral $\bm{\MF{I}(x,y,z)}$ based on an ordinary
differential equation}
\label{sect:3}

In the section we shall apply the idea of \cite{Levin82} and reduce evaluation
of the integral $\MF{I}(\vec{x})$ defined by \eqref{eq:defI} to finding one
special (bounded and slowly oscillatory) solution of an ordinary differential
equation. For the convenience of the following numerical treatment we write
$\MF{I}(\vec{x})$ as an integral over the finite interval $(0,1)$. By the change
of variable $t=\tau/(1-\tau)$ we have
\begin{equation}
  \MF{I}(\vec{x})=\int_0^1\frac{\E{\MF{\oscf_*}(\tau,\vec{x})}}{(1-\tau)^2}\,\D\tau,
\label{eq:Idef}
\end{equation}
where $\MF{\oscf_*}(\tau,\vec{x})=\MF{\oscf}(\tau/(1-\tau),\vec{x})$. It is
notable that, unlike the case considered in \cite{Levin82}, there is the
infinity of oscillations of the integrand in the neighborhood of $\tau=1$ and
the point needs a special attention.

Following \cite{Levin82}, we write the identity
\begin{equation}
  \left(\frac{\MF{\Psi}(\tau)}{(1-\tau)^2}\E{\MF{\oscf_*}(\tau,\vec{x})}\right)'_\tau=
  \frac{\E{\MF{\oscf_*}(\tau,\vec{x})}}
  {(1-\tau)^2}\left[\MF{\Psi'}(\tau)+\frac{2\MF{\Psi}(\tau)}{1-\tau}+
  \MF{\pd{\tau}\oscf_*}(\tau,\vec{x})\MF{\Psi}(\tau)\right],\label{eq:Levin0}
\end{equation}
valid for an arbitrary differentiable function $\Psi$.  Assume that we can find
$\MF{\Psi}(\tau)=\MF{\Psi}(\tau,\vec{x})$ such that the expression in the square
brackets in the last identity is equal to one. Then, obviously, we arrive at
\[
  \MF{I}(\vec{x})=\frac{\MF{\Psi}(\tau,\vec{x})}{(1-\tau)^2}\E{\MF{\oscf_*}(\tau,\vec{x})}\Bigm|_{\tau=0}^{\tau=1}.
\]

The factors of $\Psi$ in the expression in the square brackets in
\eqref{eq:Levin0} are singular at $\tau=1$, in particular, we have
\[
 \MF{\pd{\tau}\oscf_*}(\tau,\vec{x})=\frac{\MF{\sigma}(\tau,\vec{x})}{(1-\tau)^3},
\]\vskip-4mm
\noindent where\vskip-4mm
\begin{equation}
 \MF{\sigma}(\tau,\vec{x})=
 \ii x\frac{\tau(1-\tau)}{\sqrt{2\tau^2-2\tau+1}}+2y\tau+\ii z\frac{3\tau^2-2\tau+1}{\sqrt{2\tau^2-2\tau+1}}.
\label{eq:sigmadef}
\end{equation}
Therefore, assuming $|y+\ii z|\neq0$ (i.e.\ $\MF{\sigma}(1,\vec{x})\neq0$), it is
convenient to define
\[
  \MF{\Psi}(\tau,\vec{x})=(1-\tau)^3\MF{\Phi}(\tau,\vec{x}).
\]
Then, we are looking for $\MF{\Phi}(\tau,\vec{x})$ satisfying on $[0,1]$ the
following ordinary differential equation:
\begin{equation}
 (1-\tau)^3\MF{\pd{\tau}\Phi}(\tau,\vec{x})+
 \left[\MF{\sigma}(\tau,\vec{x})-(1-\tau)^2\right]\MF{\Phi}(\tau,\vec{x})=1.
\label{eq:PhiODE}
\end{equation}

Further we will show that it is possible to find a particular solution to
\eqref{eq:PhiODE} that is bounded at \mbox{$\tau=1$}. (Solutions of this type
arose in \cite{Levin82,Levin97}, where they were named slowly oscillatory, ---
in our case the particular solution is also slowly oscillatory in comparison
with the rapidly oscillatory solution to the homogeneous equation.) For the
bounded solution, taking into account the definition of $\oscf_*$, we find
\begin{equation}
  \MF{I}(\vec{x})=-\MF{\Phi}(0,\vec{x})\E{y+\ii x}.
\label{eq:IviaPhi}
\end{equation}

Let us assume that $y<0$. It is not difficult to check that the general
solution of \eqref{eq:PhiODE} can be written in the following form:
\begin{equation}
 \MF{\Phi^\textrm{gen}}(\tau,\vec{x})=c\frac{\E{\MF{\Lambda}(\tau,\vec{x})}}{\tau-1}+
 \MF{\Phi}(\tau,\vec{x}),
\label{eq:Phi_gen}
\end{equation}
where $c$ is a constant,
\begin{gather}
\MF{\Lambda}(\tau,\vec{x})=-(\tau-1)^{-2}\left(y+\ii
z\tau\sqrt{2\tau^2-2\tau+1}\right)
 -(\tau-1)^{-1}\left(2y-\ii x\sqrt{2\tau^2-2\tau+1}\right),
\label{eq:Lambdadef}\\
 \MF{\Phi}(\tau,\vec{x})=\frac{\E{\MF{\Lambda}(\tau,\vec{x})}}{\tau-1}\int_\tau^1\frac{\E{-\MF{\Lambda}(\theta,\vec{x})}}
{(\theta-1)^2}\,\D\theta.\label{eq:Phi0def}
\end{gather}

The first term in the right-hand side of \eqref{eq:Phi_gen} is highly
oscillating near $\tau=1$; besides, the amplitude of the oscillations grows when
approaching $\tau=1$ (under the assumption $y<0$). So, we fix $c=0$ in
\eqref{eq:Phi_gen} and consider the particular solution
$\MF{\Phi}(\tau,\vec{x})$.

From \eqref{eq:Lambdadef} we find
\begin{equation}
\begin{gathered}
  \MF{\Lambda}(\tau,\vec{x}) = \frac{\gamma_2}{(\tau-1)^2}+\frac{\gamma_1}{\tau-1}+
 \gamma_0+\MF{O}(\tau-1),\ \ \mbox{as}\ \ \tau\to1,\quad\mbox{where}\\
 \gamma_2=-y-\ii z,\quad \gamma_1=\ii x-2y-2\ii z,\quad
 \gamma_0=\ii x-3\ii z/2.
\end{gathered}
\label{eq:Lambdaasdef}
\end{equation}
Then we can write
\begin{equation*}
\begin{gathered}
\MF{\Phi}(\tau,\vec{x})=
\frac{\E{\MF{\Lambda}(\tau,\vec{x})}}{\tau-1}\int_\tau^1\frac{\E{-\frac{\gamma_2}{(\theta-1)^2}-\frac{\gamma_1}{\theta-1}-
 \gamma_0}}{(\theta-1)^2}\MF{\Xi}(\theta,\vec{x})\,\D\theta,
\\ \MF{\Xi}(\theta,\vec{x})=\exp\left\{-\MF{\Lambda}(\theta,\vec{x})
 +\frac{\gamma_2}{(\theta-1)^2}+\frac{\gamma_1}{\theta-1}+\gamma_0\right\}
 =1+\MF{\rho}(\theta,\vec{x}),
\end{gathered}
\end{equation*}
where $\MF{|\rho}(\theta,\vec{x})|\leq \MF{C}(\vec{x})|\theta-1|$ for $\theta\in[0,1]$ and $C$ is constant in $\theta$.
So, we have
\begin{gather}
\MF{\Phi}(\tau,\vec{x})= \frac{\E{\MF{\Lambda}(\tau,\vec{x})}}{\tau-1} \left\{
\MF{L_0}(\tau,\vec{x})+
\MF{\tilde{L}}(\tau,\vec{x})\right\},\quad\mbox{where}
\label{eq:*2.25--}\\
\MF{L_0}(\tau,\vec{x})=\int_\tau^1(\theta-1)^{-2}\E{-\frac{\gamma_2}{(\theta-1)^2}
 -\frac{\gamma_1}{\theta-1}-\gamma_0}\,\D\theta,\quad
\MF{\tilde{L}}(\tau,\vec{x})=\int_\tau^1
\MF{\rho}(\theta,\vec{x})(\theta-1)^{-2}
\E{-\frac{\gamma_2}{(\theta-1)^2}
 -\frac{\gamma_1}{\theta-1}-\gamma_0}
\,\D\theta.\notag
\end{gather}

It is not difficult to find
\begin{equation}
  \MF{L_0}(\tau,\vec{x}) = \frac{\sqrt{\pi}}{2\sqrt{\gamma_2}}\E{-\gamma_0+\frac{\gamma_1^2}{4\gamma_2}}
 \MF{\mathrm{erfc}}\left(\frac{\sqrt{\gamma_2}}{1-\tau}-\frac{\gamma_1}{2\sqrt{\gamma_2}}\right),
\label{eq:*2.5--}
\end{equation}
where we used change of variable $t=(1-\theta)^{-1}$, the equality $\gamma_2
t^2-\gamma_1
t=\gamma_2\left(t-\frac{\gamma_1}{2\gamma_2}\right)^2-\frac{\gamma_1^2}{4\gamma_2}$,
definitions of the error functions $\mathrm{erf}$ and $\mathrm{erfc}$ (see
\cite{AS}, 7.1.1, 7.1.2) and the fact that $\mathrm{erf}(\zeta)\to1$ as
$\zeta\to\infty$, $|\mathrm{arg}(\zeta)|<\pi/4$ (\cite[7.1.16]{AS}).
Additionally, for $\tau\in[0,1]$ we have
\begin{equation*}
 \MF{|\tilde{L}}(\tau,\vec{x})|\leq\MF{C}(\vec{x})(1-\tau)\int_\tau^1
 \frac{\E{\frac{y}{(\theta-1)^2}+\frac{2y}{\theta-1}}
 }{(\theta-1)^2}\D\theta=\MF{C}(\vec{x})(1-\tau)\frac{\sqrt{\pi}}{2\sqrt{|y|}}\E{-y}
 \mathrm{erfc}\left(\frac{\sqrt{|y|}\tau}{1-\tau}\right).
\end{equation*}
Then, a simple analysis based on the asymptotic expansion \cite[7.1.23]{AS} for
the $\mathrm{erfc}$ function leads us to the estimate
\begin{equation}
 \left|\frac{\E{\MF{\Lambda}(\tau,\vec{x})}}{\tau-1}\MF{\tilde{L}}(\tau,\vec{x})\right|
 =\MF{O}\bigl(\tau-1\bigr)
 \ \ \mbox{as}\ \ \tau\to1-0.
\label{eq:exp_l_est}
\end{equation}

Taking into account \eqref{eq:*2.25--}, \eqref{eq:*2.5--}, \eqref{eq:exp_l_est}
and using the expansion \cite[7.1.23]{AS} we find
\begin{gather}
\MF{\Phi}(\tau,\vec{x})=\MF{\Phi}(1,\vec{x})+
 \MF{O}(\tau-1),
\ \ \mbox{as}\ \ \tau\to1-0,\quad\mbox{where}
\notag\\
 \MF{\Phi}(1,\vec{x}):=\lim_{\tau\to1-0}\MF{\Phi}(\tau,\vec{x})=\frac{1}{2(y+\ii z)}.
 \label{eq:Phi0lim}
\end{gather}
Thus, we have also presented the boundary condition that characterizes the
bounded and slowly oscillating solution of \eqref{eq:PhiODE}.

\section{Numerical solution of (\ref{eq:PhiODE}), (\ref{eq:Phi0lim})}
\label{sect:5}

In this section we will be concerned with finding numerical solution of
\eqref{eq:PhiODE}, \eqref{eq:Phi0lim}. At this we will consider $\vec{x}$ as a
fixed parameter and for brevity in many cases we will avoid explicit denoting of
the dependence on $\vec{x}$ for $\Phi$ and some other values.

Obviously, for finding solution to \eqref{eq:PhiODE}, \eqref{eq:Phi0lim} one can
try to apply standard ODE solvers. We have tested the built-in subroutines of
the used in the present work numerical computing environment GNU Octave and the
subroutines provided by \texttt{odepkg} package to solve stiff
differential-algebraic equations. Amidst these subroutines, the most reliable
and applicable for the widest range of parameters $(x,y,z)$ was \texttt{ode5r},
which is a realization of the RADAU5 solver \cite{RADAU5}, using an implicit
Runge--Kutta method of order 5 with step size control. For application of the
subroutine, the equation \eqref{eq:PhiODE} was transformed to a system of the
form $\bm{M}\,(\Phi',\Psi')^\transp=\vec{f}(\Phi,\Psi)$ with degenerate but
constant matrix $\bm{M}$ in the left-hand side. (Here and below the symbol
$\scriptstyle\transp$ means transposition.) Namely, we use the system
\[
  \left\{
\begin{aligned}
  \Psi'={}&{}1-\left[2(1-\tau)^2+\MF{\sigma}(\tau,x,y,z)\right]\Phi,\\
  0={}&{}\Psi-(1-\tau)^3\Phi,
\end{aligned}
  \right.
\]
with the boundary conditions $\MF{\Psi}(1)=0$, $\MF{\Phi}(1)=1/(2(y+\ii z))$.
We denote by $\MF{I_\accuracy^\textrm{R5}}(x,y,z)$ the corresponding approximation
of $\MF{I}(x,y,z)$ found with \eqref{eq:IviaPhi}; here $\accuracy$ is the accuracy
demanded from the subroutine \texttt{ode5r}. Comparing to other methods
developed in the present paper for evaluation of $\MF{I}(\vec{x})$, the main
disadvantage of \texttt{ode5r} solver is its being rather slow for the
particular problem. The subject will be shortly addressed in \S\,\ref{sect:8},
but generally we will not discuss details of application of standard solvers to
\eqref{eq:PhiODE}, \eqref{eq:Phi0lim}.

In the numerical algorithm developed in the present section we seek the solution
to \eqref{eq:PhiODE} in the form of a polynomial $\MF{\Phi_M}(\tau)$ of the
order $M$. This, obviously, excludes from \eqref{eq:Phi_gen} the rapidly
oscillatory unbounded solutions of the homogeneous equation and allows us to
find an approximation of $\MF{\Phi}(\tau)$, given by \eqref{eq:Phi0def},
avoiding explicit usage of the condition \eqref{eq:Phi0lim}. For a
discretization of the differential equation and reduction it to a linear system
we use a collocation scheme demanding that $\MF{\Phi_M}(\tau)$ satisfies
\eqref{eq:PhiODE}  at the following set of Chebyshev points on the interval
$[0,1]$:
\begin{equation}
 \tau_k = (1-\cos(k\pi/M))/2,\quad k = 0, 1,..., M.
\label{eq:tausdef}
\end{equation}

We write $\MF{\Phi_M}(\tau)$ as a Lagrange interpolating polynomial so that the
values $\MF{\Phi_M}(\tau_k)\approx\MF{\Phi}(\tau_l)$ are unknowns of the linear
system. Although polynomial interpolations are usually avoided in numerical
practice, recently (see \cite{BerrutTrefethen2004,TrefethenBook,Higham}) it was
noted that the Lagrange interpolation is highly effective and numerically stable
when the polynomial is manipulated through the formulas of barycentric
interpolation and the nodes are clustered near the ends of the interval of
approximation, in particular, as for Chebyshev points.

So, an approximation of the solution to \eqref{eq:PhiODE}, \eqref{eq:Phi0lim}
is sought in the form of the barycentric Lagrange interpolating polynomial:
\begin{gather}
  \MF{\Phi_M}(\tau)=\sum_{m=0}^{\smash{M}}\MF{\ell_m}(\tau)\,\Phi_m,
\label{eq:PhiM}
\\
  \MF{\ell_m}(\tau)=\frac{\omega_m}{\tau-\tau_m}\biggm/\sum_{k=0}^M\frac{\omega_k}{\tau-\tau_k},\quad
\omega_k=\left[\prod_{l=0,\,l\neq k}^M(\tau_k-\tau_l)\right]^{-1},\quad
k,m=0,1,...,M.
\label{eq:elldef}
\end{gather}
By construction $\Phi_l=\MF{\Phi_M}(\tau_l)$, $l=0,1,...,M$. For the set of the
Chebyshev points $\{\tau_k\}$, defined by \eqref{eq:tausdef}, expressions for
$\omega_k$ are found in \cite{Salzer1972}. Cancelling factors independent of $k$
we write
\begin{equation*}
 \omega_k=(-1)^{k}\varsigma_k,\quad \varsigma_k=\begin{cases}1,&k=1\ \textrm{or}\ k=M,\\2,&\textrm{otherwise}.\end{cases}\\
\end{equation*}

Using \eqref{eq:elldef} we write the following expression
\begin{equation}
 \MF{\ell'_m}(\tau)=\frac{\omega_m}{(\tau-\tau_m)}
 \sum_{k=0}^M\frac{\omega_k}{(\tau-\tau_k)^2}\left[\sum_{k=0}^M\frac{\omega_k}{\tau-\tau_k}\right]^{-2}
 -
 \frac{\omega_m}{(\tau-\tau_m)^2}
 \left[\sum_{k=0}^M\frac{\omega_k}{\tau-\tau_k}\right]^{-1}.
\label{eq:ell'def}
\end{equation}
In view of \eqref{eq:ell'def} it is not difficult to find
\begin{gather*}
 \MF{\ell'_m}(\tau_l)=\begin{cases}\displaystyle
  \frac{\omega_m}{\omega_l(\tau_l-\tau_m)} & \mbox{for}\ m\neq l,
  \\[3mm]\displaystyle
  \sum_{k=0,\,k\neq m}^M\frac{\omega_k}{\omega_m(\tau_k-\tau_m)} & \mbox{for}\ \ m=l.
\end{cases}
\end{gather*}
This, obviously, results in the formula
\begin{equation}
  \MF{\Phi'_M}(\tau_l)=\sum_{k=0,\,k\neq l}^M\frac{\omega_k}{\omega_l(\tau_k-\tau_l)}\left[\Phi_l-\Phi_k\right].
\label{eq:PhiM'}
\end{equation}

Further we substitute \eqref{eq:PhiM}, \eqref{eq:PhiM'} into \eqref{eq:PhiODE}
and demand the equality to hold at the points $\tau_k$, $k=0,1,...,M$. In this
way we arrive at the linear system
\begin{equation}
 \bm{A}\vec{\Phi}=\vec{b},
\label{eq:linsyst}
\end{equation}
where $\vec{\Phi}=(\Phi_0,\Phi_1,...,\Phi_M)^\transp$ and
$\vec{b}=(1,1,...,1)^\transp$. To write an expression of the matrix $\bm{A}$ we
define $\mathring{\bm{A}}=[\skew4\mathring{a}_{ij}]_{i,j=0}^M$, where
$\skew4\mathring{a}_{ii}=0$ and
$\skew4\mathring{a}_{ij}=-\omega_j/(\omega_i(\tau_j-\tau_i))$. Let also
$s_i=\sum_{j=0}^M\skew4\mathring{a}_{ij}$, $i=0,1,...,M$. Then,  we have $\bm{A}
= \textrm{diag}\bigl\{(1-\vec{\tau})^3\bigr\}\cdot
\bigl(\mathring{\bm{A}}-\textrm{diag}\{s_0,s_1,...,s_M\}\bigr) +
\textrm{diag}\{\MF{\sigma}(\vec{\tau},x,y,z)-(1-\vec{\tau})^2\bigr\}$.

Finally, we solve the linear system \eqref{eq:linsyst} and in accordance with
\eqref{eq:IviaPhi} the approximation of the sought value of \eqref{eq:defI} by
the scheme based on Levin's equation and Lagrange interpolation is as follows
\begin{equation}
  \MF{I}(x,y,z) \approx \MF{I^\textrm{L--L}_M}(x,y,z):=-\Phi_0\E{y+\ii x}.
\label{eq:LevinIappr}
\end{equation}
Here is a Matlab (Octave) code of a function finding this approximation:
\begin{Verbatim}[xleftmargin=0mm,frame=leftline,fontsize=\small]
function res = Iapprox(x,y,z,M)
  sigma = @(tau,x,y,z) 2*y*tau + 1i*(x*tau.*(1-tau)+z*(3*tau.^2-2*tau+1)) ...
   ./sqrt(2*tau.^2-2*tau+1);
  tau = (1-cos(pi*(0:M)'/M))/2; w = [1; 2*ones(M-1,1); 1].*(-1).^((0:M)');
  A0 = -repmat(w',M+1,1)./repmat(w,1,M+1)./(repmat(tau',M+1,1)-repmat(tau,1,M+1));
  A0(1:(M+2):end) = 0;  A0 = A0 - diag(sum(A0,2));
  A = diag((1-tau).^3)*A0 + diag(sigma(tau,x,y,z)-(1-tau).^2);
  Phi = A\ones(M+1,1);
  res = -exp(y+1i*x) * Phi(1);
end
\end{Verbatim}

\begin{figure}[t!]\centering
  \SetLabels
  \L (0.88*-0.02) {\small $\tau$}\\
  \endSetLabels
  \leavevmode\AffixLabels{\includegraphics[width=70mm]{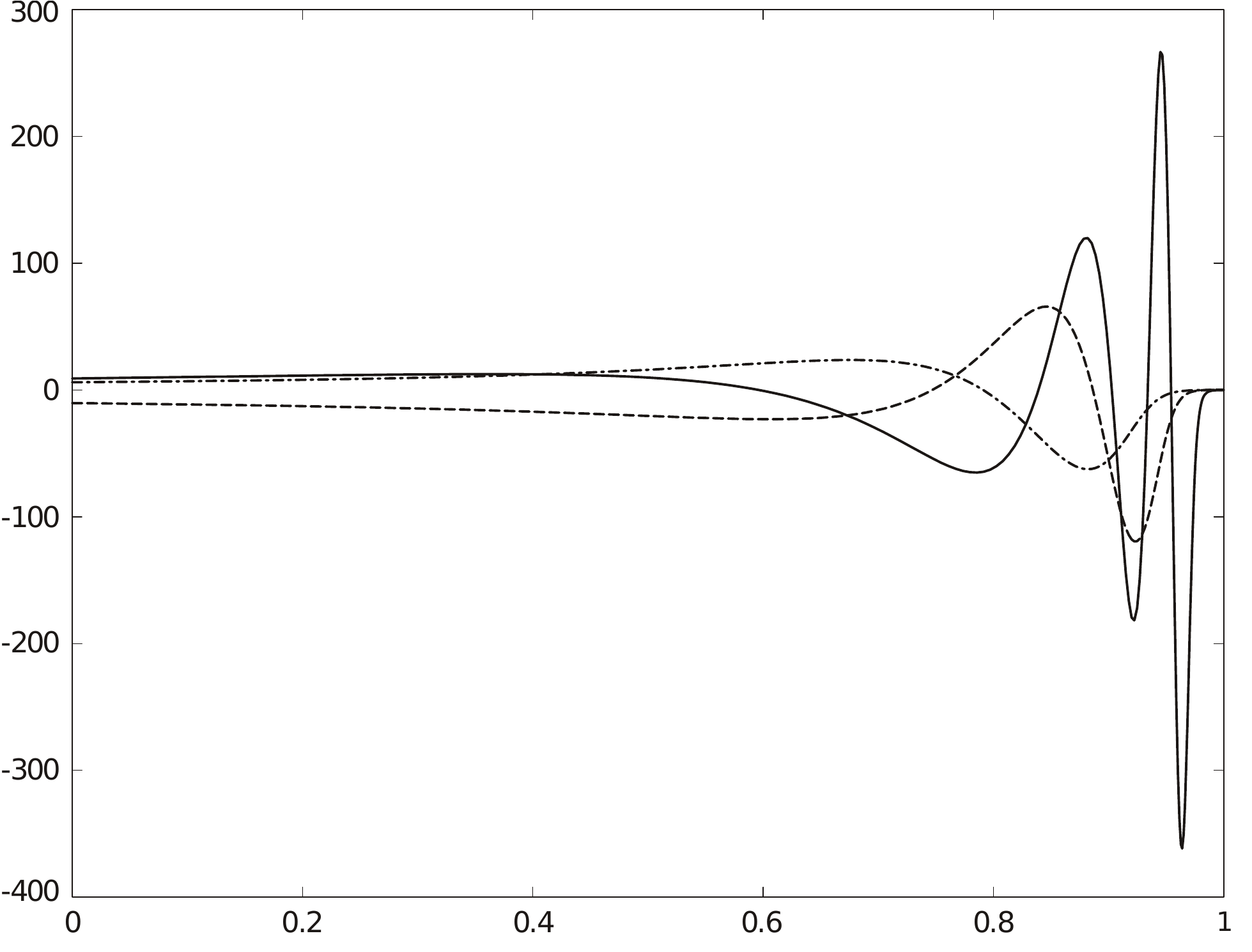}}
  \vspace{-1.5mm}
  \caption{The graph of $\MF{\Re\{\Phi}(\tau)\}$ for $x=-1$, $y=0$, $z=0.015$ (solid line), 
  $z=0.03$ (dashed line), $z=0.045$ (dash-dotted line).}
  \label{fig:phi}
\end{figure}

To improve the above scheme we can seek the solution to \eqref{eq:PhiODE}, \eqref{eq:Phi0lim} in
the form
\begin{equation}
  \MF{\Phi}(\tau)\approx\MF{\skew{5}\Hat\phi}(\tau)+\MF{\Phi_M}(\tau),
\label{eq:PhiM_Tm_Mod}
\end{equation}
where $\MF{\Phi_M}(\tau)$ is defined by \eqref{eq:PhiM} and $\skew{5}\Hat\phi$
is a function that is intended to absorb some irregularities in behaviour of
$\Phi$ (like sharp peaks). As shown in fig.~\ref{fig:phi}, in particular, the
irregular behaviour is characteristic when $y=0$, $z>0$ and the point $(x,y,z)$
approaches $(x,0,0)$ (track behind the source located on the free surface). The peaks concentrate near the critical value
$\tau_*=\bigl(2z-x+\sqrt{x^2-8z^2}\bigr)/(6z-2x)$, where
$\Im\{\MF{\sigma(\tau_*,\vec{x})}\}=0$.

Since $\tau_*$ is close to $1$ when $z$ is small, we will use the representation
\eqref{eq:*2.25--} and in view of \eqref{eq:Lambdaasdef} we define
\[
  \MF{\skew{5}\Hat\phi}(\tau)=\frac{\E{\frac{\gamma_2}{(\tau-1)^2}+\frac{\gamma_1}{\tau-1}+
 \gamma_0}}{\tau-1}\MF{L_0}(\tau,\vec{x}).
\]
It is convenient to rewrite the last formula in terms of Faddeeva function
$\MF{\mathrm{w}}(\eta)=\E{-\eta^2}\mathrm{erfc}(-\ii \eta)$
\cite{FaddeyevaTerentiev1961}. By a simple algebra we find
\begin{equation}
\MF{\skew{5}\Hat\phi}(\tau)=\frac{\sqrt{\pi}}{2(\tau-1)\sqrt{\gamma_2}}
\MF{\mathrm{w}}\left(\frac{\ii\sqrt{\gamma_2}}{1-\tau}-\frac{\ii\gamma_1}{2\sqrt{\gamma_2}}\right).
\label{eq:hatphidef}
\end{equation}
The last function can be effectively evaluated using known algorithms and
packages (e.g.\ \cite{FaddeevaPackage}). 

By using \eqref{eq:PhiM_Tm_Mod} we arrive at the linear system
\eqref{eq:linsyst} with another right-hand side:
\begin{equation}
\vec{b}^*=(1,1,...,1)^\transp-\vec{\ell},
\label{eq:bstdef}
\end{equation}
where
$\vec{\ell}=\bigl(\MF{\Hat{L}}(\tau_0),\MF{\Hat{L}}(\tau_1),...,\MF{\Hat{L}}(\tau_{M})\bigr)^\transp$,
\begin{multline}
\!\!\!\MF{\Hat{L}}(\tau)=\mathcal{L}\MF{\skew{5}\Hat\phi}(\tau)=(1-\tau)^3\MF{\skew{5}\Hat\phi'}(\tau)+
\bigl[\MF{\sigma}(\tau,x,y,z)-(1-\tau)^2\bigr]\MF{\skew{5}\Hat\phi}(\tau)\\{}=
 1+\frac{\sqrt{\pi}}{2\sqrt{y+\ii z}(\tau-1)}
 \left\{\frac{\bigl(3\tau^2-2\tau+1\bigr)z+\tau(1-\tau)x}{\sqrt{2\tau^2-2\tau+1}}
 -2\tau z+(\tau-1)x \right\}
\MF{\mathrm{w}}\left(\frac{\tau\sqrt{y+\ii z}}{\tau-1}-\frac{\ii x}{2\sqrt{y+\ii z}}\right).
\label{eq:l0phidef}
\end{multline}

Further we shall denote by $\MF{\hat{I}^\textrm{L--L}_M}(x,y,z)$ the
approximation of $\MF{I}(x,y,z)$ obtained with the help of
\eqref{eq:PhiM_Tm_Mod} and \eqref{eq:hatphidef}. Practically,
$\MF{\hat{I}^\textrm{L--L}_M}(x,y,z)$ has advantages comparing with
$\MF{I^\textrm{L--L}_M}(x,y,z)$ when $\tau_*$ is close to $1$.

It is also useful to have a way to control errors of the approximation.
After solving the system \eqref{eq:linsyst}, by using \eqref{eq:PhiM} and
defining $\MF{\Phi'_M}(\tau)$ with the help of \eqref{eq:ell'def} we can find
the residual $\MF{r}(\tau)$ of the equation \eqref{eq:PhiODE}
($r=\mathcal{L}\Phi_M-1$). Analogously, if we solve the system
$\bm{A}\vec{\Phi}=\vec{b}^*$, where $\vec{b}^*$ is given by \eqref{eq:bstdef},
then the residual of the equation \eqref{eq:PhiODE} is defined by
$r=\mathcal{L}\bigl(\Phi_M+\skew{5}\Hat\phi\bigr)-1$. At this, the residual of
the solution $R=\Phi_M-\Phi$ (or $R=\Phi_M+\skew{5}\Hat\phi-\Phi$) satisfies the
equation
\begin{equation}
 \mathcal{L}\mkern1muR=r.
\label{eq:ee*}
\end{equation}

Since $\MF{R}(\tau)$ is the slowly oscillating, bounded solution to
\eqref{eq:ee*}, following arguments of \S\,\ref{sect:3}, we find
\begin{equation}
  \MF{R}(\tau)=\frac{\E{\MF{\Lambda}(\tau)}}{\tau-1}
 \int_\tau^1\frac{\MF{r}(\theta)\E{-\MF{\Lambda}(\theta)}}{(\theta-1)^2}\,\D\theta,
\label{eq:eeR*}
\end{equation}
where $\Lambda$ is defined in \eqref{eq:Lambdadef}. The representation
\eqref{eq:eeR*} allows us to obtain a simple estimates of
$\bigl|I-I^\textrm{L--L}_M\bigr|$ and $\bigl|I-\hat{I}^\textrm{L--L}_M\bigr|$ if
$y<0$. Namely, fixing $\tau=0$ in \eqref{eq:eeR*} and taking into account
\eqref{eq:IviaPhi}, \eqref{eq:Lambdadef}, \eqref{eq:LevinIappr} we find
\begin{equation}
 \bigl|I-I^\textrm{L--L}_M\bigr|\leq
 \E{2y}\left|\int_0^1\frac{\MF{r}(\theta)\E{-\MF{\Lambda}(\theta,\vec{x})}}{(1-\theta)^2}\,\D\theta\right|
 \leq\E{2y}\max_{\theta\in[0,1]}\{|r(\theta)|\}
 \int_0^1\frac{\E{\frac{y}{(\theta-1)^2}+\frac{2y}{\theta-1}}}{(\theta-1)^2}\D\theta\\{}=
  \frac{\sqrt{\pi}}{2\sqrt{|y|}}\E{y}\max_{\theta\in[0,1]}\{|r(\theta)|\}.
\label{eq:ee**}
\end{equation}
This estimate also holds for $\bigl|I-\hat{I}^\textrm{L--L}_M\bigr|$. In
practice, it is convenient to compute $r_k=\MF{r}(\check{\tau}_k)$ at
$\check{M}=M-1$ Chebyshev points
\begin{equation*}
\check\tau_k =
\frac12\left(1-\cos\left(\frac{(k+1/2)\pi}{\check{M}+1}\right)\right),\quad
k=0,1,...,\check{M}.
\end{equation*}
At this $\check\tau_k \in (\tau_k,\tau_{k+1})$ and taking into account the
nature of $r(\tau)$ which oscillates between zero values at $\tau_k$, we can
assume that $\max_{k=0,1,...,M}\{|r_k|\}$ approximates the term
$\max_{\theta\in[0,1]}\{|r(\theta)|\}$ in the right-hand side of
\eqref{eq:ee**}.

Unfortunately, the estimate \eqref{eq:ee**} is poor for small $|y|$. In this
case we can use the following approach. We approximate $R(\tau)$ by
\begin{equation}
  \MF{R_{\check{M}}}(\tau)=\sum_{m=0}^{\check{M}}\MF{\skew4\check\ell_m}(\tau)\,R_m,
\label{eq:RM}
\end{equation}
where $\MF{\skew4\check\ell_m}(\tau)$ is computed analogously to
\eqref{eq:elldef}, for $\check{\tau}_k$ and $\check{\omega}_k$ defined as
follows (see e.g.\ \cite{BerrutTrefethen2004}):
\begin{equation*}
 \check\omega_k=(-1)^{k}\sin\bigl((k+1/2)\pi/\bigl(\check{M}+1\bigr)\bigr),
 \quad k=0,1,...,\check{M}.
\end{equation*}

Substituting \eqref{eq:RM} and the corresponding representation of
$R'_{\smash{\check{M}}\vphantom{M}}$ (see \eqref{eq:PhiM'}) into \eqref{eq:ee*} we
demand the equality to hold at the points $\check{\tau}_k$,
$k=0,1,...,\check{M}$. In this way we obtain the linear system
\begin{equation}
 \bm{B}\vec{R}=\vec{r},
\label{eq:linsysterr}
\end{equation}
where $\vec{R}=(R_0,R_1,...,R_{\check{M}})^\transp$,
$\vec{r}=(r_0,r_1,...,r_{\check{M}})^\transp$. Defining
$\mathring{\bm{B}}=\bigl[\mathring{b}_{ij}\bigr]_{i,j=0}^{\check{M}}$,
$\mathring{b}_{ii}=0$,
$\mathring{b}_{ij}=-\check{\omega}_j/(\check{\omega}_i(\check{\tau}_j-\check{\tau}_i))$
and $\check{s}_i=\sum_{j=0}^{\check{M}}\mathring{b}_{ij}$, we have $\bm{B} =
\textrm{diag}\bigl\{(1-\check{\vec{\tau}})^3\bigr\}\cdot
\bigl(\mathring{\bm{B}}-\textrm{diag}\{\check{s}_0,\check{s}_1,...,\check{s}_M\}\bigr)
+
\textrm{diag}\{\MF{\sigma}(\check{\vec{\tau}},x,y,z)-(1-\check{\vec{\tau}})^2\bigr\}$.

We solve the system \eqref{eq:linsysterr} and assume that
$|R(0)|\leq\|\vec{R}\|_\infty=\max\{|R_k|:k=0,1,...,\check{M}\}$ (assumption
$|R(0)|\leq|R_0|$ would be more accurate but less reliable choice). Finally, we
can estimate $\bigl|I-I^\textrm{L--L}_M\bigr|$ and
$\bigl|I-\hat{I}^\textrm{L--L}_M\bigr|$ by the value
\begin{equation}
 \MF{\varepsilon^\textrm{L--L}_M}(x,y,z)=\E{y}\min\left\{
 \|\vec{R}\|_\infty,
 2^{-1}(\pi/|y|)^{1/2}\|\vec{r}\|_\infty\right\}.
\label{eq:eeeps}
\end{equation}
Our numerical experiments (see \S\,\ref{sect:8}) show that the estimate
$\varepsilon_M^\textrm{L--L}$ is rather reliable. However, evaluation of
$\varepsilon_M^\textrm{L--L}$ is time-consuming, so in \S\,\ref{sect:8} we also
discuss ways to avoid the computations.

\section[Evaluation of $\MF{I}(x,y,z)$ based on Clenshaw--Curtis quadrature and
steepest descent method]{Evaluation of $\bm{\MF{I}(x,y,z)}$ based on Clenshaw--Curtis quadrature and
steepest descent method}
\label{sect:7}

Another scheme, developed in the present work to compute $\MF{I}(x,y,z)$, is
based on application of Clen\-shaw--Curtis quadrature where, as a preliminary, we
use the steepest descent method. (It is notable that the latter method could be
used to improve the scheme of \S\,\ref{sect:3}, \ref{sect:5}.)

For the function $\MF{\oscf}(t,x,y,z)$ involved in the definition of
$\MF{I}(x,y,z)$ by \eqref{eq:defI}, we have
\begin{equation}
\MF{\oscf}(t,x,y,z) =
 (y+\ii z)t^2+\ii x t + y + \MF{O}\bigl(t^{-1}\bigr)\ \ \mbox{as}\ \ t\to\infty.
\label{eq:gasympt}
\end{equation}
Therefore, Jordan's lemma allows rotation $t \rightarrow \E{\ii\theta}t$ of the
integration contour in \eqref{eq:defI} if $\Re\bigl[(y+\ii
z)\E{2\ii\vartheta}\bigr]<0$ for $\vartheta\in(0,\theta)$. The value of $\theta$
corresponding to the steepest descent for sufficiently large $t$ can be found by
minimizing $\Re\bigl[(y+\ii z)\E{2\ii\theta}\bigr]$. This gives
\begin{equation}
\begin{gathered}
 \cos2\theta=\frac{-y}{\sqrt{y^2+z^2}},\ \ \sin2\theta=\frac{z}{\sqrt{y^2+z^2}},\ \ 
 \Re\bigl[(y+\ii z)\E{2\ii\theta}\bigr]=-\sqrt{y^2+z^2},
\\
 \cos\theta=\sqrt{\frac{1+|y|/\sqrt{y^2+z^2}}{2}},
  \ \ \sin\theta=\textrm{sign}(z)\sqrt{\frac{1-|y|/\sqrt{y^2+z^2}}{2}}.
\end{gathered}
\label{eq:2theta}
\end{equation}

We note that $\textrm{sign}(z)\,\theta\in[0,\frac{\pi}{4}]$. Therefore, for
$z<0$, $\Re\bigl[\ii x \E{\ii\theta}\bigr]<0$ and after the rotation the second
term in the right-hand side of \eqref{eq:gasympt} provides the decaying
exponential factor $\exp\{-xt\sin\theta\}$ in the integrand as $t\to\infty$.
Hence, for $z\leq0$ we will use the following expression:
\begin{equation}
  \MF{I}(x,y,z)=\E{\ii\theta}\int_0^{\infty}\E{\MF{\oscf}(\E{\ii\theta}t,x,y,z)}\,\D{}t.
\label{eq:Irot1}
\end{equation}

The case $z>0$ needs more attention, because then $\Re\bigl[\ii x
\E{\ii\theta}\bigr]>0$ and the second term in the right-hand side of
\eqref{eq:gasympt} prevails in some finite interval of $t\geq0$. So, the
definition \eqref{eq:Irot1} would not be useful for numerical computations due
to the presence of very large (in absolute magnitude) values of the integrand.
Therefore, we split the interval of integration into two parts and write
\begin{equation}
  \MF{I}(x,y,z)=\int_0^{t_*}\E{\MF{\oscf}(t,x,y,z)}\,\D{}t+
  \E{\ii\theta}\int_0^{\infty}\E{\MF{\oscf}(\MF{\gamma}(t),x,y,z)}\,\D{}t,
\label{eq:Irot2}
\end{equation}
where $\MF{\gamma}(t)=t_*+t\E{\ii\theta}$. We need the real part of
$\MF{\oscf}(\MF{\gamma}(t),x,y,z)$ to be negative for $t\geq0$. Looking at
\eqref{eq:gasympt}, to find $t_*$ we demand that for $t\geq0$
\begin{equation}
  |x|\MF{\Im}(\MF{\gamma}(t))\leq z\MF{\Im}\bigl([\MF{\gamma}(t)]^2\bigr)
  +|y|\bigl[1+\MF{\Re}\bigl([\MF{\gamma}(t)]^2\bigr)\bigr].
\label{eq:ineq4tst}
\end{equation}
By using \eqref{eq:2theta} we rewrite the right-hand side of the latter inequality as
follows:
\[
 2tt_*\left(|y|\cos\theta+z\sin\theta\right)+t^2\sqrt{y^2+z^2}+|y|\left(1+t_*^2\right).
\]
Thus, it is easy to find that the inequality \eqref{eq:ineq4tst} holds for
$t\geq0$ if $|x|\MF{\Im}(\MF{\gamma}(t))\leq
2tt_*\left(|y|\cos\theta+z\sin\theta\right)$. Thus, in our computations we fix
\begin{equation}
  t_*=\frac{|x|\sin\theta}{2(|y|\cos\theta+z\sin\theta)}.
\label{eq:t*def}
\end{equation}

To find approximation of the integrals in \eqref{eq:Irot1}, \eqref{eq:Irot2} we
will use the Clenshaw--Curtis quadrature \cite{ClCu}. First we give an outline
of the method emphasizing the possibility of effective application of the
quadrature using FFT. Consider the integral
\begin{equation}
 \int_{-1}^1\MF{f}(t)\,\D{}t
\label{eq:intf-1_1}
\end{equation}
for some function $f:[-1,1]\mapsto\RR$. Expand the function $f$ into the series
of Chebyshev polynomials:
\begin{equation}
\MF{f}(t)=\sum_{k=0}^\infty\epsilon_k c_k \MF{T_k}(t),
\label{eq:chext}
\end{equation}
where $\epsilon_0=1/2$, $\epsilon_k=1$, $k\geq1$,
$\MF{T_k}(t)=\cos(k\arccos(t))$. Then, it can be found that
\begin{equation}
 \int_{-1}^1\MF{f}(t)\,\D{}t=c_0+2\sum_{k=1}^\infty\frac{c_{2k}}{1-(2k)^2}=
 d_0+2\sum_{k=1}^\infty\frac{d_{k}}{1-(2k)^2},
\label{eq:f_in_ds}
\end{equation}
where we define $d_k=c_{2k}$, $k=0,1,2,...$ In view of the definition
$\MF{T_k}(t)$, the change $t=\cos(\tau)$ in \eqref{eq:chext} and the
orthogonality of the cosines in $\MF{L_2}(0,\pi)$ lead us to
\[
 d_k=\frac{2}{\pi}\int_0^\pi\MF{f}(\cos\tau)\cos(2k\tau)\,\D{}\tau.
\]

Let us now fix a sufficiently large, even integer $N$ and use the trapezoidal
rule to approximate the latter integral:
\begin{equation}
 d_k\approx\delta_k=\frac{2}{N}\left[
 \frac{\MF{f}(1)}{2}+\frac{\MF{f}(-1)}{2}+
 \sum_{n=1}^{N-1}\MF{f}(\cos(n\pi/N))\cos(2\pi nk/N)
 \right].
\label{eq:dapprox}
\end{equation}
It is easy to note that \eqref{eq:dapprox} expresses $\delta_k$ through the
discrete Fourier transform that maps one vector of length $N$ (say, $\vec{u}$)
to another one (say, $\vec{U}$) via $U_{k,N}=\sum_{n=0}^{N-1}u_n\E{-\ii 2\pi
kn/N}$. We write $\vec{U}=\mathcal{F}(\vec{u})$ or, equivalently, in the matrix
notation $\vec{U}=\bm{F}\vec{u}$, where
$\bm{F}=[\phi_{\alpha\beta}]_{\alpha,\beta=1}^N$,
$\phi_{\alpha\beta}=\E{-2\pi\ii(\alpha-1)(\beta-1)/N}$. Then, we have
\begin{equation}
 \vec{\delta}=(\delta_0,\delta_1,...,\delta_{N-1})^\transp=\frac{2}{N}\Re\{\bm{F}\}\vec{f}_*,
\label{eq:deltaasF}
\end{equation}
where
$\vec{f}_*=\left([\MF{f}(1)+\MF{f}(-1)]/2,\MF{f}(t_{1,N}),...,\MF{f}(t_{N-1,N})\right)^\transp$,
$t_{k,N}=\cos k\pi/N$ are extrema of $\MF{T_N}(t)$.

We note that $\delta_{\frac{N}{2}+\ell}=\delta_{\frac{N}{2}-\ell}$,
$\ell=1,2,...,N/2-1$, and define the vector
$\vec{\kappa}=(1,\kappa_1,...,\kappa_{N-1})^\transp$, such that
$\kappa_n=1/\bigl(1-4n^2\bigr)$ for $n=1,2,...,N/2$ and
$\kappa_n=1/\bigl(1-4(N-n)^2\bigr)$ for $n=N/2+1,N/2+2,...,N-1$. Then, in view
of \eqref{eq:f_in_ds}, \eqref{eq:dapprox} and \eqref{eq:deltaasF} we can write
\[
  \int_{-1}^1\MF{f}(t)\,\D{}t\approx \delta_0+2\sum_{k=1}^{N/2}\frac{\delta_k}{1-4k^2}=
  \vec{\kappa}^\transp\vec{\delta}=\frac{2}{N}\vec{\kappa}^\transp\Re\{\bm{F}\}\vec{f}_*=
  \frac{2}{N}\bigl(\Re\{\bm{F}\}\vec{\kappa}\bigr)^\transp\vec{f}_*,
\]
where we use the fact that $\Re\{\bm{F}\}$ is symmetric. Finally, we find
\begin{equation}
  \int_{-1}^1\MF{f}(t)\,\D{}t\approx
  \vec{v}^\transp\vec{f}_*=\vec{w}^\transp\vec{f},
\label{eq:C-Cdef}
\end{equation}
where
$\vec{v}=\frac{2}{N}\Re\{\bm{F}\}\vec{\kappa}=\frac{2}{N}\Re\{\mathcal{F}(\vec{\kappa})\}$
and the second equality is the usual form of a quadrature with weights
$\vec{w}=(2^{-1}v_0,v_1,...,v_{N-1},2^{-1}v_0)^\transp$ and values of the
integrand
$\vec{f}=\bigl(\MF{f}(t_{0,N}),\MF{f}(t_{1,N}),...,\MF{f}(t_{N,N})\bigr)^\transp$.
The weights $\vec{v}$ and $\vec{w}$ can be effectively computed using FFT.
(Applicability of FFT for calculation of weights of Clenshaw--Curtis quadrature
was first pointed out in \cite{Gentleman1972I}.) So, the quadrature shares the
usual excellent numerical stability of the FFT (see e.g.\ \cite{KanekoLiu1970})
and the low cost of implementation.

Further we transform the integrals in \eqref{eq:Irot1}, \eqref{eq:Irot2}, where
$t_*$ is defined by \eqref{eq:t*def}, to integrals over $[-1,1]$ of the type
\eqref{eq:intf-1_1} by the following changes of variables
\begin{equation}
\begin{gathered}
 \int_0^{t_*}\MF{q}(t)\,\D{}t=\frac{t_*}{2}\int_{-1}^1\MF{q}\left(\frac{t_*(t+1)}{2}\right)\,\D{}t=
 \int_{-1}^1\MF{f}(t)\,\D{}t,
\\
 \int_0^\infty\MF{q}(t)\,\D{}t=2\int_{-1}^1(1-t)^{-2}\MF{q}\left(\frac{t+1}{1-t}\right)\,\D{}t=
 \int_{-1}^1\MF{f}(t)\,\D{}t,
\end{gathered}
\label{eq:chvar4C-C}
\end{equation}
Consider any one of the integrals arising in \eqref{eq:Irot1}, \eqref{eq:Irot2}
after the change of variables \eqref{eq:chvar4C-C}. We introduce the notation
\begin{equation}
  \MF{\mathring{F}_N}(f)=\vec{w}^\transp\vec{f}
\label{eq:Fndef}
\end{equation}
for the approximation obtained through \eqref{eq:C-Cdef}. It is easy to note
that for $y^2+z^2\neq0$ the function $f$ in question is absolutely continuous on
$[-1,1]$ along with all derivatives $f^{(k)}$, $k=1,2,...$. Thus, by Theorem~5.1
\cite{Trefethen} for any integer $k>0$ and all sufficiently large $N$
\begin{equation}
 \rho_N:=\left|\MF{\mathring{F}_N}(f)-\int_{-1}^1\MF{f}(t)\,\D{}t\right|\leq C N^{-k},
\label{eq:estimr}
\end{equation}
where the constant $C$ depends on $f$ and $k$.

We will apply the approximation \eqref{eq:Fndef} in a `brute force' manner,
increasing $N$ until we are sure that $\rho_{N}\leq\accuracy$ for the given
accuracy $\accuracy$. The reasons, why one can expect to obtain a relatively
effective scheme in this way, are the guaranteed by \eqref{eq:estimr}
convergence of the sequence $\mathring{F}_N=\MF{\mathring{F}_N}(f)$ to its limit
value as $N\to\infty$; the insensitivity of Clenshaw--Curtis scheme to rounding
errors in the evaluation of the integrand; low ($N\log(N)$) cost of
implementation and simplicity of the integrand. As we shall see in
\S\,\ref{sect:8}, the suggested numerical approach works rather well even
sufficiently close to the line $y=0$, $z=0$.

Since the Clenshaw--Curtis quadrature is naturally nested, in order not to waste
already computed integrand values when computing the sequence $\{\mathring{F}_N\}$,
it is convenient to choose the next step value of $N$ to be a multiple of the
previous one. Hence we will consider $F_\ell=\mathring{F}_{2^\ell N_0}$,
$\ell=0,1,2,...$, for some initial even integer $N_0$ (in computations,
presented in \S\,\ref{sect:8}, $N_0=2$).

A number of estimates for $\rho_N$ have been proposed~---~see
\cite{Gentleman1972I,Gentleman1972II} and references therein. We shall use the simple
error estimate, demanding convergence of $F_\ell$ as a Cauchy sequence --- so that the
computation is stopped at the step $\ell=\ell_*$ if the last $\#$ members of the sequence
$\{F_\ell\}$ are sufficiently close to each other:
$\max\left\{\left|F_{\ell_*-i+1}-F_{\ell_*-j+1}\right|:1\leq i\leq\#,\, 1\leq j\leq\#\right\}
 \leq\accuracy$.
It is known that for rapidly convergent processes such estimate can be rough,
essentially overestimating the error. However, as it is emphasized in
\cite{Gentleman1972I} the considerable experience of implementation of the
Clenshaw--Curtis scheme confirms that the simple estimate is very realistic. The
conclusion agrees with the observations made in the present work.

In our computations, presented in the next section, $\#=3$ and we also
introduce the reserve coefficient $c_r=10$ for the difference of $F_{\ell_*}$
and $F_{\ell_*-1}$, so that the termination condition is as follows:
\begin{equation}
 \max\left\{c_r\left|F_{\ell_*}-F_{\ell_*-1}\right|,\left|F_{\ell_*}-F_{\ell_*-2}\right|,
 \left|F_{\ell_*-1}-F_{\ell_*-2}\right|\right\}\leq\accuracy.
\label{eq:ce_used}
\end{equation}

Further we will denote by $\MF{I^\textrm{C--C}_\accuracy}(x,y,z)$ the
approximation of $\MF{I}(x,y,z)$ obtained with the help of \eqref{eq:Irot1},
\eqref{eq:Irot2}, \eqref{eq:C-Cdef}, \eqref{eq:chvar4C-C}, \eqref{eq:Fndef} and
\eqref{eq:ce_used}. We also denote by
$\MF{\mathcal{N}^\textrm{C--C}_\accuracy}(x,y,z)$ the total number of
evaluations of the functions under integral signs in \eqref{eq:Irot1} or
\eqref{eq:Irot2}, needed to satisfy \eqref{eq:ce_used}. For $z\leq0$, when we
use the expression \eqref{eq:Irot1},
$\MF{\mathcal{N}^\textrm{C--C}_\accuracy}(x,y,z)=2^{\ell_*}N_0+1$; for $z>0$,
$\MF{\mathcal{N}^\textrm{C--C}_\accuracy}(x,y,z)=\left(2^{\ell'_*}+2^{\ell''_*}\right)N_0+2$,
where $\ell'_*$ ($\ell''_*$) is such that \eqref{eq:ce_used} holds for the first
(second) integral in the right-hand side of \eqref{eq:Irot2}.

\section{Numerical experiments}
\label{sect:8}

In this section we present results of numerical testing of the approaches
developed in \S\S\,\ref{sect:3}, \ref{sect:5} and \ref{sect:7}. We study
properties of the approximations of $\MF{I}(x,y,z)$, with the particular
attention to their accuracy and efficacy. Calculations are performed in the
numerical computing environment GNU Octave and double precision (64-bit)
arithmetics (the machine epsilon $\epsilon$ is about $2.22\cdot10^{-16}$).

\begin{figure}[t!]\vspace{2mm}
\begin{center}
 \SetLabels
  \L (0.76*0.51) $x$\\
  \L (0.76*-0.01) $x$\\
  \L (0.0*0.99) $z$\\
  \L (0.0*0.47) $z$\\
  \L (1.0*0.89) $\log_{10}(\mathcal{N}^\textrm{C--C}_\accuracy)$\\
  \endSetLabels
  \mbox{\kern-1mm}\AffixLabels{\includegraphics[width=110mm]{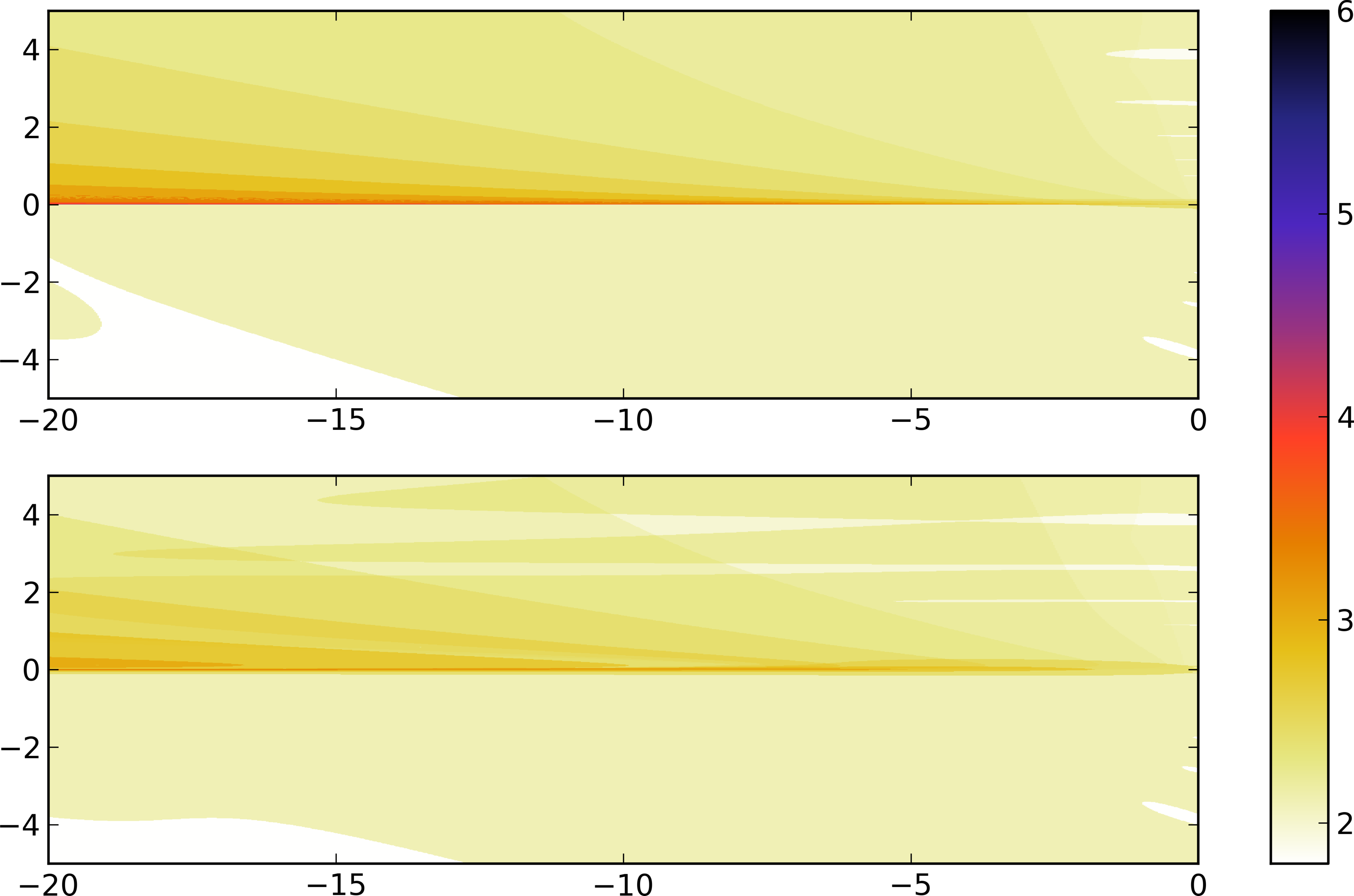}}
\end{center}\vspace{-4mm}
\caption{Values of
$\MF{\log_{10}}\bigl(\MF{\mathcal{N}^\textrm{C--C}_\accuracy}(x,y,z)\bigr)$ for
$\accuracy=10^{-6}$ and $y=0$ (upper), $y=-0.1$ (lower).}
\label{fig:num1}
\end{figure}

\begin{figure}[t!]
\vspace{3mm}
\begin{center}
 \SetLabels
  \L (0.76*0.51) $x$\\
  \L (0.76*-0.01) $x$\\
  \L (0.02*0.99) $z$\\
  \L (0.02*0.47) $z$\\
  \L (1.0*0.89) $\log_{10}(\mathcal{N}^\textrm{C--C}_\accuracy)$\\
  \endSetLabels
  \mbox{\kern-2.5mm}\AffixLabels{\includegraphics[width=112mm]{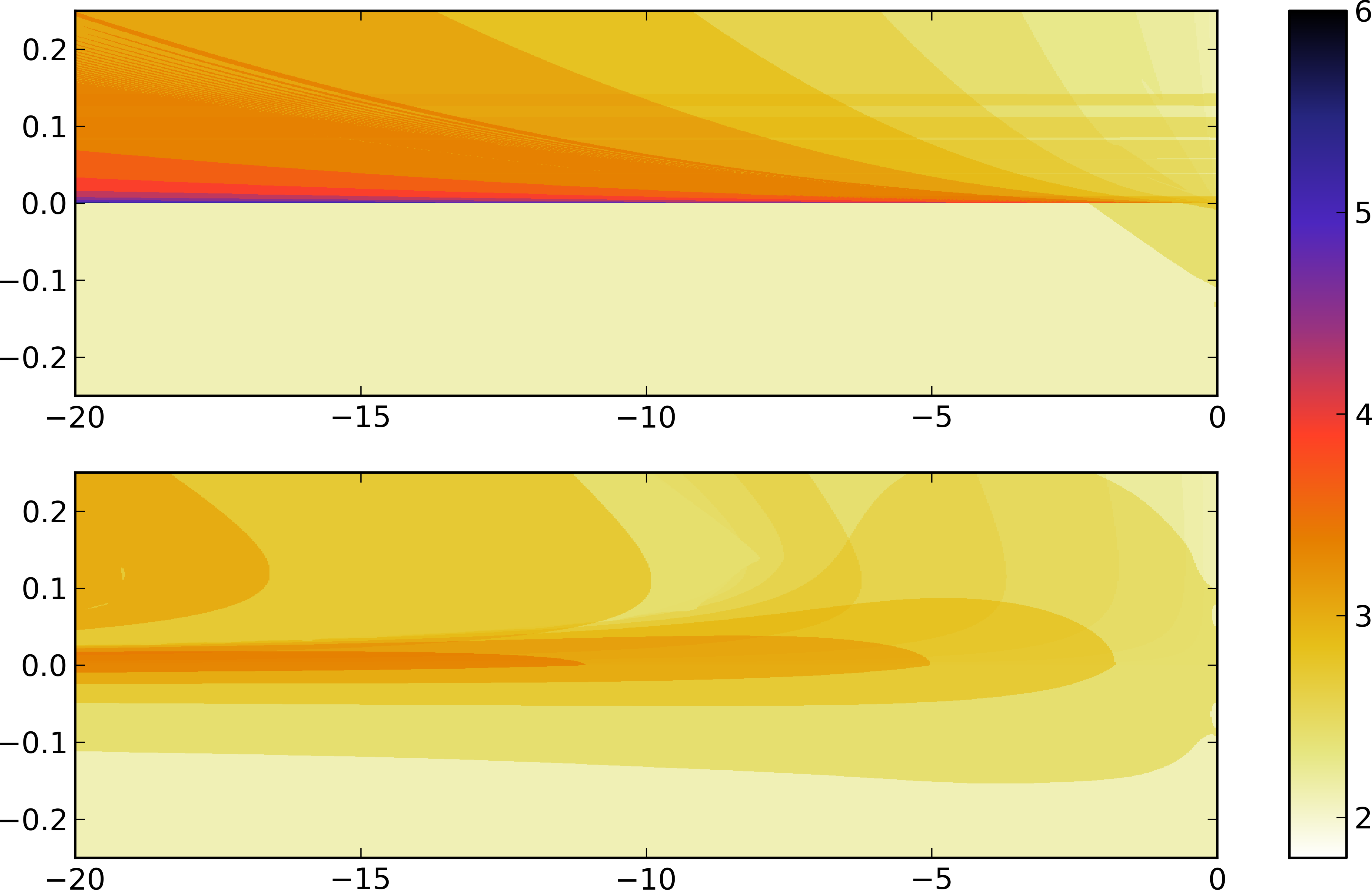}}
\end{center}\vspace{-4mm}
\caption{Values of
$\MF{\log_{10}}\bigl(\MF{\mathcal{N}^\textrm{C--C}_\accuracy}(x,y,z)\bigr)$ for
$\accuracy=10^{-6}$ and $y=0$ (upper), $y=-0.1$ (lower).}
\label{fig:num2}
\end{figure}

\begin{figure}[t!]\vspace{2mm}
\begin{center}
 \SetLabels
  \L (0.76*0.51) $x$\\
  \L (0.76*-0.01) $x$\\
  \L (0.02*0.99) $z$\\
  \L (0.02*0.47) $z$\\
  \L (1.0*0.89) $\log_{10}(\mathcal{N}^\textrm{C--C}_\accuracy)$\\
  \endSetLabels
  \mbox{\kern-2mm}\AffixLabels{\includegraphics[width=107mm]{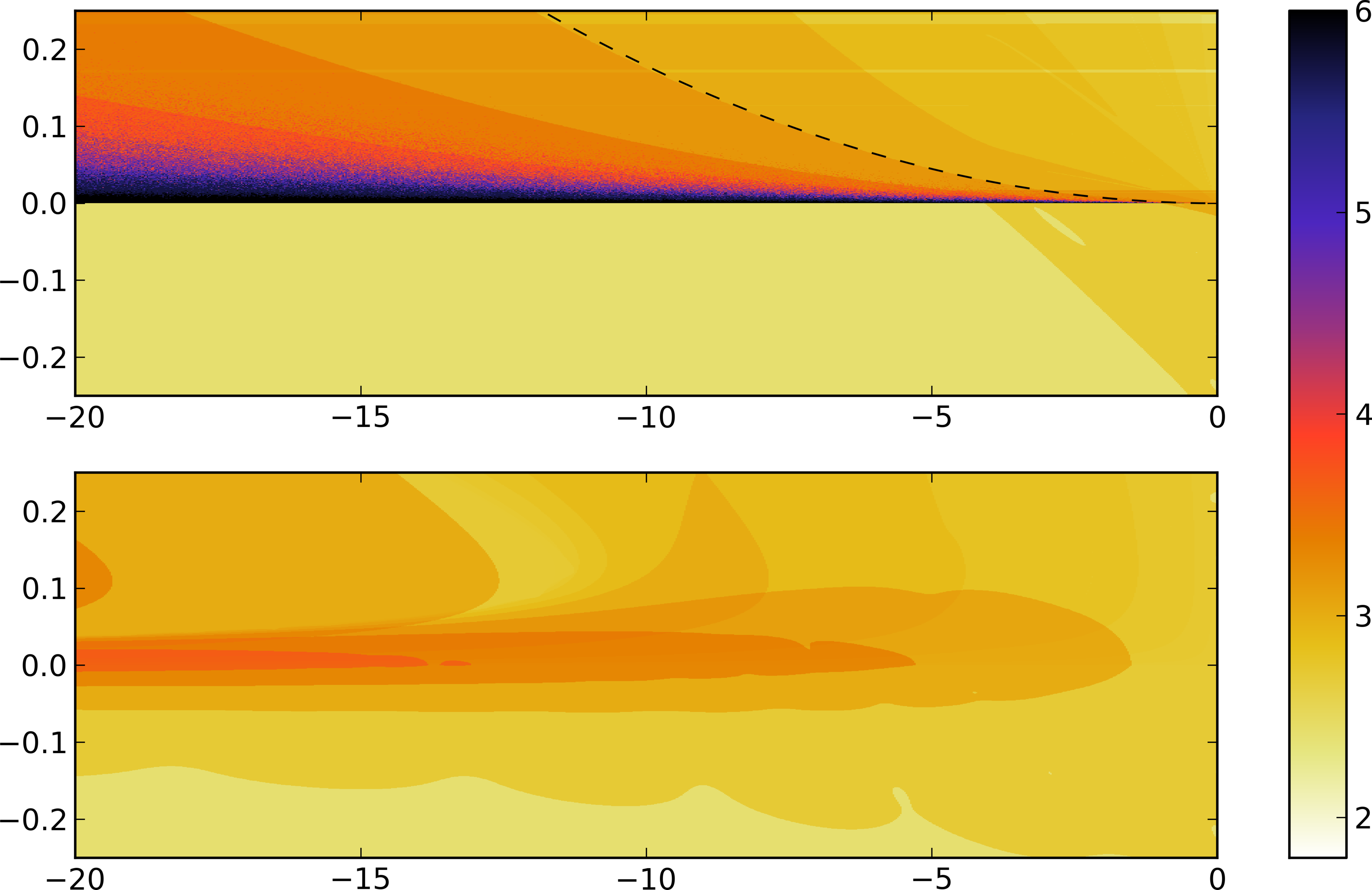}}
\end{center}\vspace{-4mm}
\caption{Values of
$\MF{\log_{10}}\bigl(\MF{\mathcal{N}^\textrm{C--C}_\accuracy}(x,y,z)\bigr)$ for
$\accuracy=10^{-12}$ and $y=0$ (upper), $y=-0.1$ (lower).}
\label{fig:num3}
\end{figure}

We start presenting computations of the quantity
$\MF{\mathcal{N}^\textrm{C--C}_\accuracy}(x,y,z)$.
Figures~\ref{fig:num1}--\ref{fig:num3} show in a semi-logarithmic scale the
values of $\MF{\mathcal{N}^\textrm{C--C}_\accuracy}(x,0,z)$ and
$\MF{\mathcal{N}^\textrm{C--C}_\accuracy}(x,-0.1,z)$ computed on the grid $(x,z)
\in\vec{L}(800,[-20,0])\times\vec{L}(800,[-5,5])$ (fig.~\ref{fig:num1}) and
$(x,z) \in\vec{L}(800,[-20,0])\times\vec{L}(800,[-0.25,0.25])$
(fig.~\ref{fig:num2}, \ref{fig:num3}). Here and below $\vec{L}(n,\gamma)$ is the
set consisting of $n$ linearly spaced in the interval $\gamma$ values (including
interval's end-points). The demanded accuracy of computation $\accuracy=10^{-6}$
in fig.~\ref{fig:num1}, \ref{fig:num2}, and $\accuracy=10^{-12}$ in
fig.~\ref{fig:num3}.

\begin{figure}[t!]
\vspace{4mm}
\begin{center}
 \SetLabels
  \L (0.02*0.99) $z$\\
  \endSetLabels
  \mbox{\kern-2mm}\AffixLabels{\includegraphics[width=80mm]{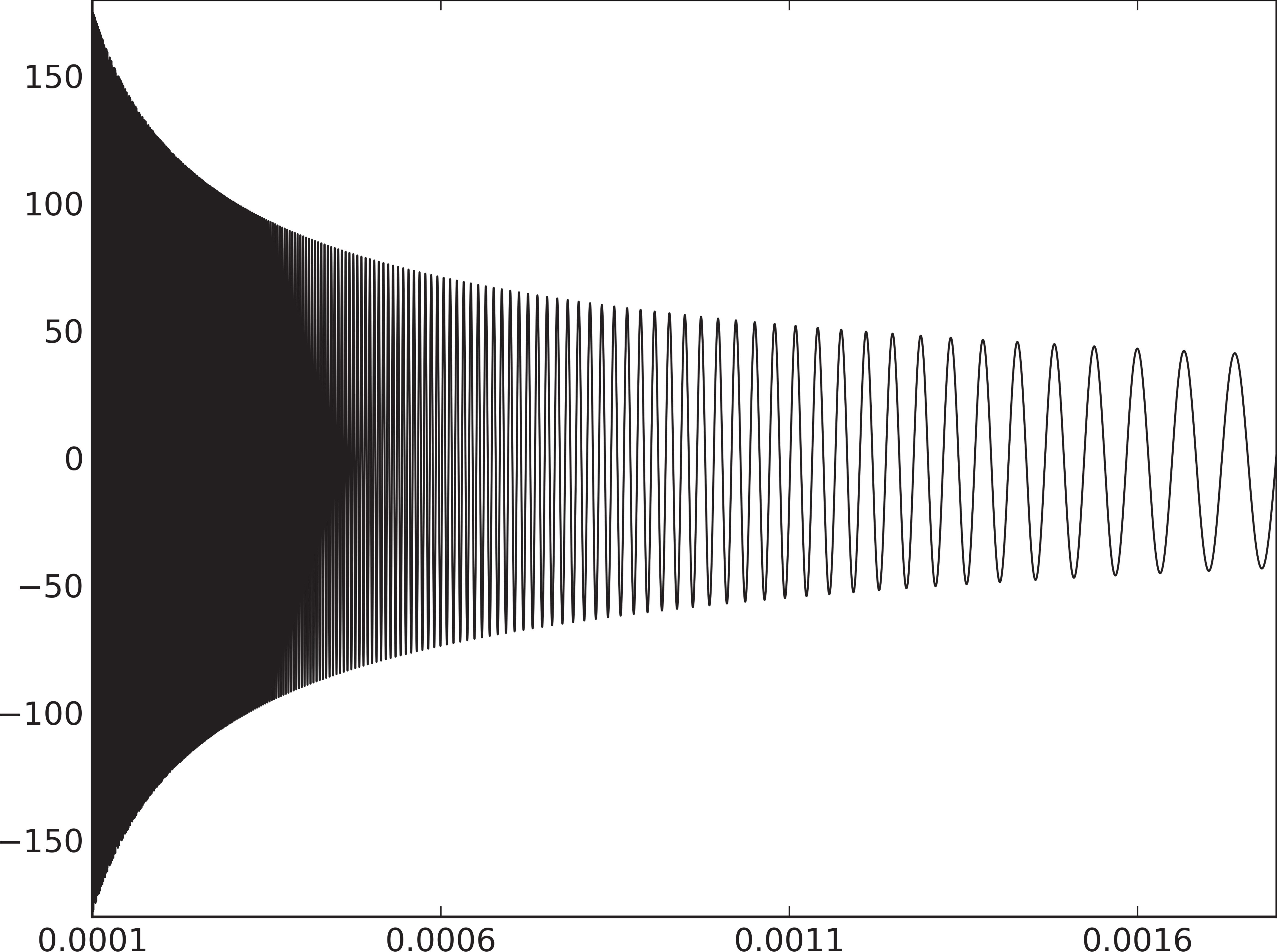}}
\end{center}\vspace{-4mm}
\caption{Dependence of $\Im\{\MF{I}(x,y,z)\}$ on $z$ for fixed $x=-1$, $y=0$.}
\label{fig:num5}
\end{figure}

In our computations we set the limit of $2^{19}+1$ evaluations of integrand for
each integral in \eqref{eq:Irot1}, \eqref{eq:Irot2}. The points, where the
accuracy is not achieved with the workload (the inequality \eqref{eq:ce_used} is
not satisfied), can be easily recognized in the upper part of
fig.~\ref{fig:num3} (the black area located above and close to the $x$-axis).
From the arguments of \S\,\ref{sect:7} (see, in particular, \eqref{eq:gasympt})
it is easy to note that the numerical integration of \eqref{eq:Irot1},
\eqref{eq:Irot2} become more difficult with the increase of either of the
parameters $D=x^2/(4\sqrt{y^2+z^2})$ or $\theta\in[-\pi/4,\pi/4]$ defined by
\eqref{eq:2theta}. This happens due to weakening of decay at infinity of the
oscillating integrands and naturally leads to accumulation of round-off errors
(this effect is more evident for the higher demanded accuracy
$\accuracy=10^{-12}$). We also note the irregularity of $\MF{I}(x,y,z)$ for
large $D$: fig.~\ref{fig:num5} shows the highly oscillating, with increasing
amplitude, behaviour of $\MF{\Im\{I}(x,0,z)\}$ as $z\to+0$. (Shown in the figure
values are computed as $\Im\bigl\{\MF{I^\textrm{C--C}_{10^{-6}}}(x,y,z)\bigr\}$.)

However, we should note that the evaluation of
$\MF{I^\textrm{C--C}_\accuracy}(x,y,z)$ is reasonably quick and effective up to
rather large values of $D$. For example, the dashed line in fig.~\ref{fig:num3}
corresponds to $D=100$ (notably, this is the worst case $\theta=\pi/4$ and the
accuracy $\accuracy=10^{-12}$ is quite close to $\epsilon$). The case of large
$D$ has been investigated extensively in \cite{Newman3}, where various
expansions of $\MF{I_\infty}(x,y,z;y_0)$ (see \eqref{eq:defI_inf}) are given and
algorithms for numerical evaluation of the function are presented. It is claimed
in \cite{Newman3} that the algorithms are applicable and effective for $D>40$,
so they are complementary to our approaches.

Of course, to substantiate the arguments above we should check reliability of
the test \eqref{eq:ce_used}. First, for this purpose we compare
$\MF{I^\textrm{C--C}_{\accuracy_i}}(x,y,z)$ for $\accuracy_1=10^{-6}$ and
$\accuracy_2=10^{-12}$ to verify that
\begin{equation}
 \left|\MF{I^\textrm{C--C}_{\accuracy_1}}(x,y,z)-\MF{I^\textrm{C--C}_{\accuracy_2}}(x,y,z)\right|\leq\accuracy_1.
\label{eq:Icc-Icc}
\end{equation}
The functions $\MF{I^\textrm{C--C}_{\accuracy_i}}(x,y,z)$ are evaluated for $y=0$,
$-0.1$, $-0.25$, $-0.5$ on the grid $(x,z)\in
\vec{L}(800,[-20,0])\times\vec{L}(800,[-0.25,0.25])$. Amidst the checked
$2.56\cdot10^6$ pairs $I^\textrm{C--C}_{\accuracy_1}$, $I^\textrm{C--C}_{\accuracy_2}$
(minus $6854$ where one or two of the values cannot be computed due to the limit
on the number of integrand evaluations) the condition \eqref{eq:Icc-Icc} is not
satisfied at $8$ points $(x,z)$ only. These points correspond to $y=-0.1$,
belong to a vicinity of radius $0.05$ of $(x,z)=(-15.92,0.032)$ (so $D$ is quite
large, about $603.47$) and the maximum found value of the expression in the
left-hand side of \eqref{eq:Icc-Icc} is approximately equal to
$1.85\cdot10^{-6}$. Hence, these computations confirm that the criterion of
termination \eqref{eq:ce_used} is sufficiently reliable.

We can also compare $\MF{I^\textrm{C--C}_{\accuracy}}(x,y,z)$ and the
approximation $\MF{I^\textrm{R5}_{\accuracy}}(x,y,z)$ introduced in
\S\,\ref{sect:5} through an application of the RADAU5 solver \cite{RADAU5} to
find the solution to \eqref{eq:PhiODE}, \eqref{eq:Phi0lim}. It is notable that
computation of $\MF{I^\textrm{R5}_{\accuracy}}(x,y,z)$ with \texttt{ode5r} GNU
Octave function is rather stable, failing only for large values of $D$ (e.g.\ it
is possible to compute $\MF{I_\accuracy^\textrm{R5}}(-1,0,0.008)$ with an
accuracy $\accuracy$ better than $10^{-10}$). Unfortunately, it is also rather
slow, for example, achieving the accuracy $\accuracy=10^{-6}$ by
$\MF{I^\textrm{R5}_\accuracy}(-0.1,0,0.01)$ is about $260$ slower than the
computation of $\MF{I^\textrm{C--C}_{10^{-6}}}(-0.1,0,0.01)$. Hence, comparison
of $\MF{I^\textrm{C--C}_{\accuracy}}(x,y,z)$ and
$\MF{I^\textrm{R5}_{\accuracy}}(x,y,z)$ for large sets of $(x,y,z)$ would be too
time-consuming. However, to give an example, we evaluate the difference
$\MF{\bigl|I^\textrm{R5}_{\accuracy_1}}(x,y,z)-\MF{I^\textrm{C--C}_{\accuracy_2}}(x,y,z)\bigr|$
for $x$ varying from $-10$ to $0$ with step $0.01$ and for fixed $y=0$, $-0.1$;
$z=0.1$, $0.2$. Since the build-in error estimate in \texttt{ode5r}
systematically underestimates
$\MF{|I}(x,y,z)-\MF{I^\textrm{R5}_\accuracy}(x,y,z)|$, we choose
$\accuracy_1=10^{-10}$, $\accuracy_2=10^{-6}$  and find the maximal value of the
difference
$\MF{\bigl|I^\textrm{R5}_{\accuracy_1}}(x,y,z)-\MF{I^\textrm{C--C}_{\accuracy_2}}(x,y,z)\bigr|$
to be about $1.58\cdot10^{-7}$.

\begin{figure}[t!]\vspace{2mm}
\begin{center}
 \SetLabels
  \L (0.83*-0.02) $x$\\
  \L (-0.02*0.5) $z$\\
  \endSetLabels
  \mbox{\kern0mm}\AffixLabels{\includegraphics[width=158mm]{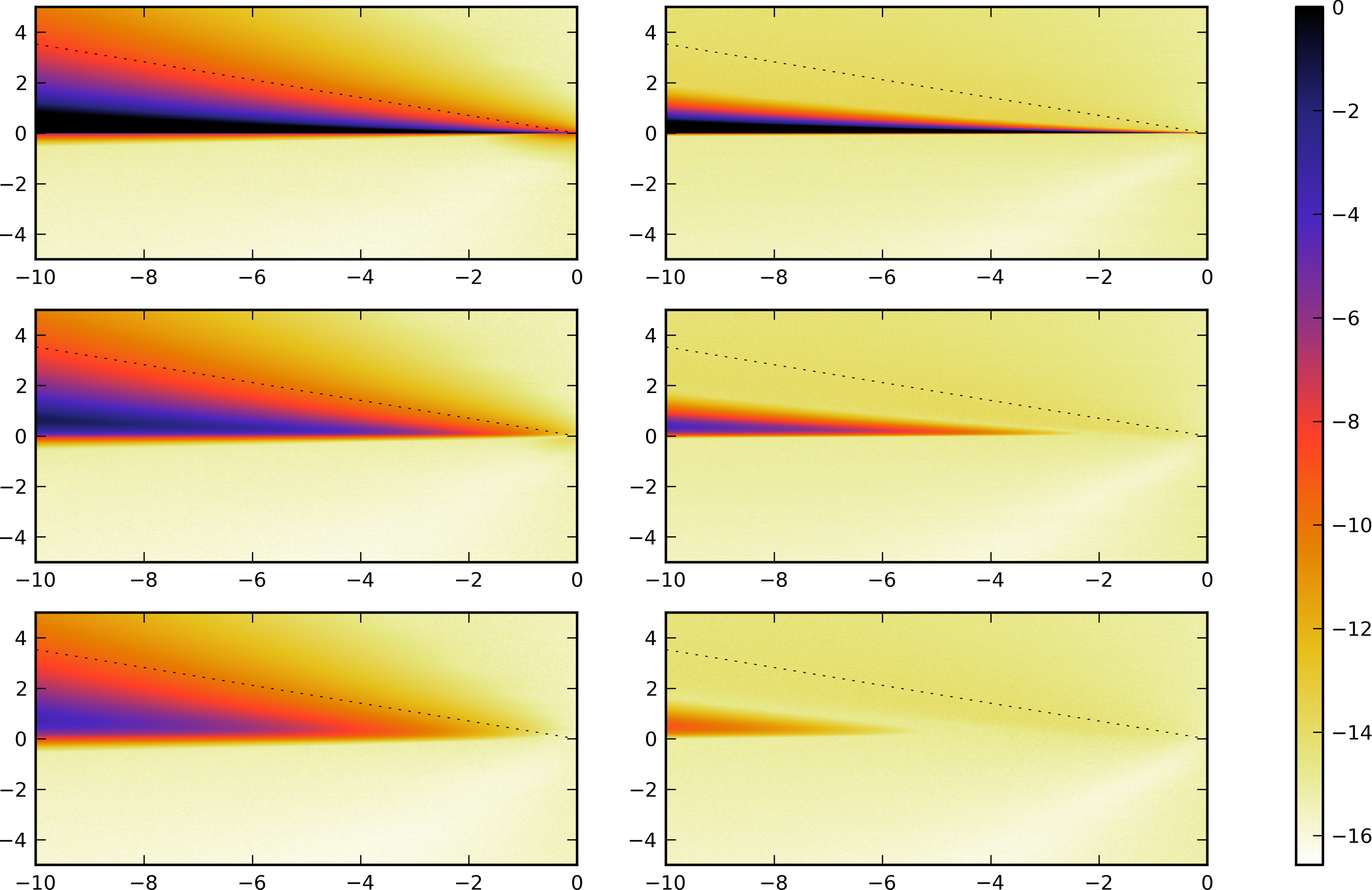}}
\end{center}\vspace{-4mm}
\caption{Computation of
$\min\left\{0,\log_{10}(\varepsilon_M^\textrm{L--L})\right\}$ as a function of
$(x,z)$ for three fixed values $y=0$, $-0.1$, $-0.25$ (from top to bottom) and
for $M=50$ (left column), $100$ (right column). The dashed line is locus of
equation $x=-2\sqrt{2}z$ (Kelvin angle).}
\label{fig:num6}
\end{figure}

\begin{figure}[t!]\vspace{2mm}
\begin{center}
 \SetLabels
  \L (0.83*-0.02) $x$\\
  \L (-0.02*0.5) $z$\\
  \endSetLabels
  \mbox{\kern0mm}\AffixLabels{\includegraphics[width=158mm]{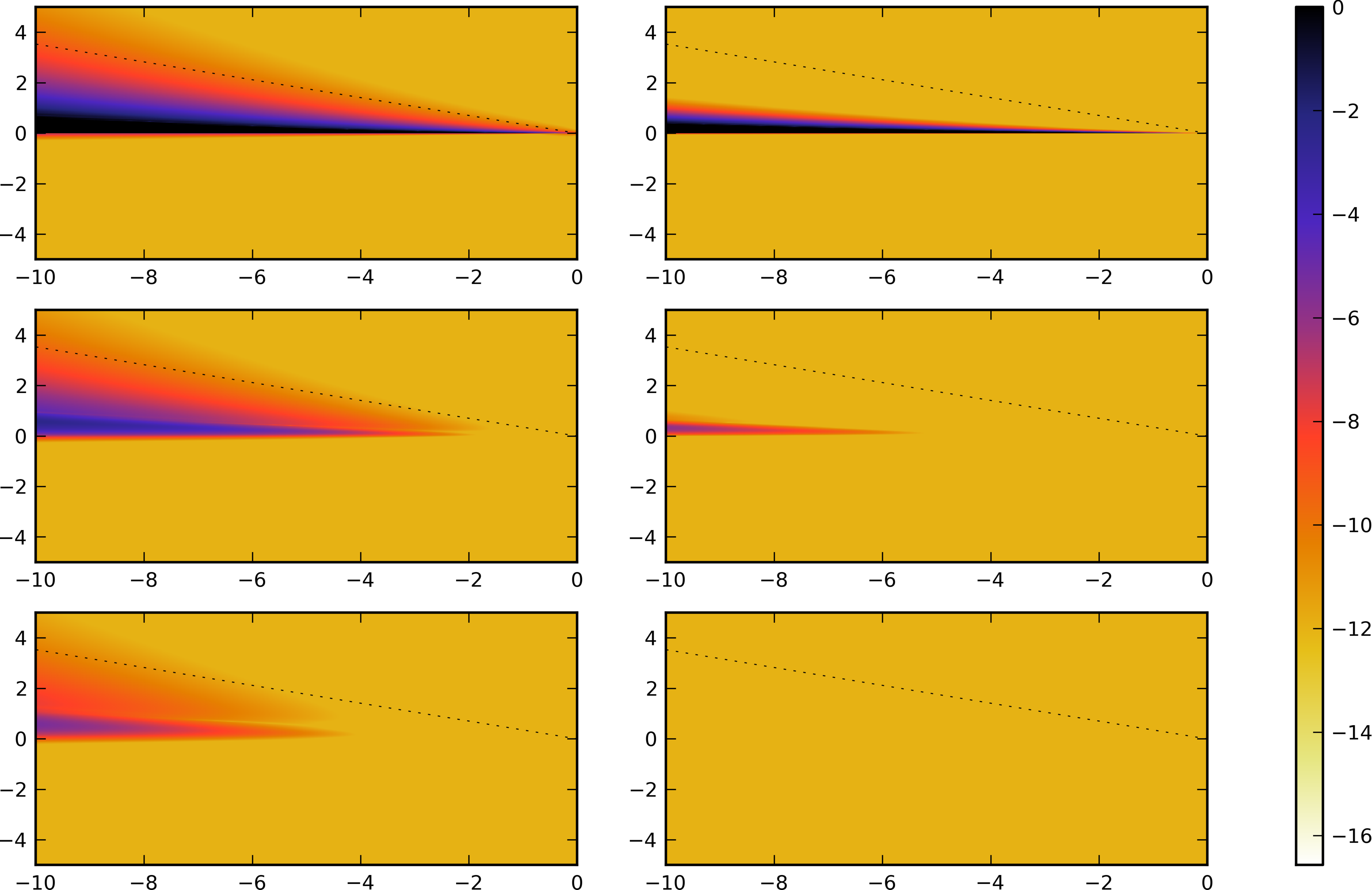}}
\end{center}\vspace{-4mm}
\caption{Computation of
$\max\Bigl\{-12,\min\Bigl\{0,\log_{10}\Bigl(\bigl|\MF{I^\textrm{C--C}_{10^{-12}}}(x,y,z)-\MF{\hat{I}^\textrm{L--L}_{M}}(x,y,z)\bigr|\Bigr)\Bigr\}\Bigr\}$
as a function of $(x,z)$ for three fixed values $y=0$, $-0.1$, $-0.25$ (from top
to bottom) and for $M=50$ (left column), $100$ (right column). The dashed line
is locus of equation $x=-2\sqrt{2}z$ (Kelvin angle).}
\label{fig:num7}
\end{figure}

Now we compare the two schemes developed in \S\S\,\ref{sect:5}, \ref{sect:7} and check the
error estimate \eqref{eq:eeeps}. For these purposes we numerically test the inequality
\begin{equation}
 \left|\MF{I^\textrm{C--C}_{\accuracy}}(x,y,z)-\MF{\hat{I}^\textrm{L--L}_{M}}(x,y,z)\right|\leq
 \max\left\{\MF{\varepsilon^\textrm{L--L}_M}(x,y,z),\accuracy\right\}
\label{eq:CCLGchk}
\end{equation}
for $y=0$, $-0.1$, $-0.25$, $-0.5$, $M=50$, $100$, and $\accuracy=10^{-12}$. The
functions are evaluated on the grid $(x,z)\in
\vec{L}(400,[-10,0])\times\vec{L}(400,[-5,5])$. Amidst the checked
$1.28\cdot10^6$ pairs $I^\textrm{C--C}_{\accuracy}$, $\hat{I}^\textrm{L--L}_{M}$
corresponding to different values of $(x,y,z)$ and $M$ we found $3145$ such that
\eqref{eq:CCLGchk} does not hold. All these cases correspond to $y=0$ and occur
near the $x$-axis, for $z>0$ and $D>49$. Besides, all the cases of violation of
\eqref{eq:CCLGchk} occur for large values of the error estimate
$\varepsilon^\textrm{L--L}_M$ (bigger than $1.15$ in our computations). Thus,
since we have already got evidence that the test \eqref{eq:ce_used} is
trustworthy, we can claim that the estimate $\varepsilon^\textrm{L--L}_M$ is
fairly adequate and properly overestimates the error
$\MF{\bigl|I}(x,y,z)-\MF{\hat{I}^\textrm{L--L}_{M}}(x,y,z)\bigr|$ providing the
accuracy of approximation is demanded to be sufficiently high (as it is typical
in practice).

Fig.~\ref{fig:num6} shows dependence of
$\MF{\varepsilon^\textrm{L--L}_M}(x,y,z)$ on $(x,z)$ for fixed values $y=0$,
$-0.1$, $-0.25$, $M=50$, $100$, in a semi-logarithmic scale. It is interesting
to compare the picture with fig.~\ref{fig:num7} showing
\[\max\left\{-12,\min\left\{0,\log_{10}\left(\left|\MF{I^\textrm{C--C}_{10^{-12}}}(x,y,z)-\MF{\hat{I}^\textrm{L--L}_{M}}(x,y,z)\right|\right)\right\}\right\}.\]
The values in fig.~\ref{fig:num6} and \ref{fig:num7} are computed on the grid
$(x,z)\in\vec{L}(400,[-10,0])\times\vec{L}(400,[-5,5])$. It is important to note
that unlike the standard Levin's collocation, our scheme, based on barycentric
Lagrange interpolation, demonstrates numerical stability, its quality increases
with growth of $M$ in the representation~\eqref{eq:PhiM_Tm_Mod} and it can be
applied for $M$'s equal to hundreds and thousands. For example,
$\MF{\hat{I}^\textrm{L--L}_{1000}}(-1,0,0.005)$ provides approximation of $I$
with an accuracy better than $10^{-12}$. Besides, we have rather satisfactory
results for subsets of the domain between the $x$-axis and the dashed line in
fig.~\ref{fig:num6}, \ref{fig:num7}~---~in the presence of critical points of
the oscillating integrand in \eqref{eq:defI}. The situation of the presence of
critical points is usually avoided in applications of Levin's collocation
scheme.

We present some information on the speed of the considered algorithms. To be
definitive we state that the computer in use has a 3.4\;GHz Intel Core i5 CPU
and 8\;Gb of RAM. Then, the evaluation time of
$\MF{\hat{I}^\textrm{L--L}_M}(x,y,z)$ is about $1.25\cdot10^{-3}$\;sec.\ for
$M=50$, $2.5\cdot10^{-3}$\;sec.\ for $M=100$, $9\cdot10^{-3}$\;sec.\ for
$M=200$. Time of evaluation of $\MF{I^\textrm{C--C}_\accuracy}(x,y,z)$ depends on
the number of integrands evaluations in \eqref{eq:Irot1}, \eqref{eq:Irot2},
needed to achieve the accuracy $\accuracy$; the time is about $10^{-2}$ when
$\log_{10}\mathcal{N}^\textrm{C--C}_\accuracy\approx 4.2$ and about
$2.5\cdot10^{-3}$ when $\log_{10}\mathcal{N}^\textrm{C--C}_\accuracy\approx 2.4$
(see Figures~\ref{fig:num1}--\ref{fig:num3}).

The experience of computations shows that the procedures
$\MF{\hat{I}^\textrm{L--L}_M}(x,y,z)$ and
$\MF{I^\textrm{C--C}_\accuracy}(x,y,z)$ spend comparable time to provide an
approximation of $\MF{I}(x,y,z)$ with the given accuracy $\accuracy$. It can be
observed that the algorithm based on Levin's ODE and barycentric Lagrange
interpolation is considerably faster than the counterpart when the value of $M$,
needed to achieve the given accuracy, is small. The situation of small $M$ is,
in particular, typical for large $|y|$; for instance, we need $M=20$ for $x=-1$,
$y=-1$, $z=0.1$, $\accuracy=10^{-12}$ and $\MF{\hat{I}^\textrm{L--L}_M}(x,y,z)$
is more than $3$ times faster. However, unlike the algorithm based on
Clenshaw--Curtis quadrature~---~designed to provide the result with the demanded
accuracy~---~in the algorithm based on Levin's ODE and Lagrange interpolation
the number $M$ guaranteeing the accuracy $\accuracy$ is not known a-priori.
However, fig.~\ref{fig:num6} and \ref{fig:num7} show that the approximation
error by $\MF{\hat{I}^\textrm{L--L}_M}(x,y,z)$ is quite regular and predictable.
Consider a fixed $\accuracy$ and test $\MF{I^\textrm{L--L}_M}(x,y,z)$ against
$\MF{I^\textrm{C--C}_{\accuracy/2}}(x,y,z)$ for a grid of $(x,y,z)$. For each
point $(x,y,z)$ of the grid we find minimal $\MF{M_\accuracy}(x,y,z)$
guaranteeing that
$\MF{\bigl|\MF{I^\textrm{C--C}_{\accuracy/2}}}(x,y,z)-\MF{\hat{I}^\textrm{L--L}_M}(x,y,z)\bigr|\leq\accuracy/2$.
(This should be a rather trivial but time-consuming operation.) Suppose that a
simple majorizing function
$\MF{\mathcal{M}_\accuracy}(x,y,z)\geq\MF{M_\accuracy}(x,y,z)$ is proposed, then
we could use
$\MF{\hat{I}^\textrm{L--L}_{\MF{\mathcal{M}_\accuracy}(x,y,z)}}(x,y,z)$ to
provide the result with the given accuracy $\accuracy$. However, such study is
beyond of our scope in the present paper.

Along with the convenience to provide the approximation of $\MF{I}(x,y,z)$ with
the given accuracy, another advantage of the algorithm based on the
Clenshaw--Curtis quadrature is its ability to work rather well (naturally,
becoming slower) for large $D$ and $\theta$ close to $\pi/4$. As
fig.~\ref{fig:num2}, \ref{fig:num3} show, for reasonable accuracy, say,
$\accuracy=10^{-6}{-}10^{-7}$~---~not allowing the limitations of the
double-float arithmetics to manifest themselves~---~the values of $D$ can be
very large (on the standards of \cite{Newman3}). To give an example, we stay
within the limit of $2^{19}+1$ evaluations of each integrand in \eqref{eq:Irot2}
when computing $\MF{I^\textrm{C--C}_{10^{-7}}}\bigl(x,0,10^{-6}\bigr)$ for $x$
varying from $-1$ to $0$ with the step $10^{-5}$. (Note that
$\MF{D}(-1,0,10^{-6})$ is as large as $2.5\cdot10^{5}$.)

\begin{table}[t!]
  {\footnotesize\centering
  \begin{tabular}{|c|c|c|c|c|}
\hline
       \diagbox{$z$}{$y$}       & $-0.5$ & $-0.1$ & $-0.01$ & $0$ \\
\hline
     $0.5$ & $-0.3132089735$ & $-0.4347821474$ &
    $-0.4093149760$ & $-0.4039184710$ \\
\hline
     $0.1$ & $-0.4288349681$ & $-1.0716691716$ &
    $-2.1157417380$ & $-2.5160949098$ \\
\hline
     $0.01$ & $-0.4349760923$ & $-0.9188289512$ &  
  $-0.7896492217$ & $ 3.6856412628$ \\
\hline
  \end{tabular}\par}
  \caption{Benchmark values of $\MF{I_\infty}(x,y,z)$ defined by \eqref{eq:defI_inf} for $x=-1$ and a variety
of $(y,z)$.}
  \label{tab:1}
\end{table}

Finally, we can conclude that both procedures
$\MF{I^\textrm{C--C}_\accuracy}(x,y,z)$ and
$\MF{\hat{I}^\textrm{L--L}_M}(x,y,z)$ are sufficiently fast and reliable. These
algorithms are complementary to those of \cite{Newman3}, developed for large
$D$. The algorithm based on Clenshaw--Curtis quadrature is rather universal,
provides results with the given accuracy and is applicable for very large values
of $D$. The algorithm based on Levin's ODE and the barycentric Lagrange
interpolation has advantages when the number $M$ needed for achieving a given
accuracy $\accuracy$ is sufficiently small.


\section[Computation of derivatives of $\MF{I}(x,y,z)$]{Computation of derivatives of $\bm{\MF{I}(x,y,z)}$}
\label{sect:der}

In this section we consider application of the suggested methods to computation
of derivatives of $\MF{I}(x,y,z)$ and show that the generalization is rather
straightforward. Let us introduce a vector
$\vec{\ell}=(\ell_1,\ell_2,\ell_3)^\transp$ and consider
\begin{equation}
 \MF{J}(\vec{\ell},\vec{x}):=\MF{\nabla I}(x,y,z)\cdot\vec{\ell}=
 \int_0^\infty \MF{\oscf}(t,\ell_1,\ell_2,\ell_3) \E{\MF{\oscf}(t,x,y,z)}\,\D{}t,
\label{eq:Jldef}
\end{equation}
where the function $\oscf$ is defined by \eqref{eq:defI}.

First, we note that in view of smoothness and polynomial behaviour at infinity
of the function $\MF{\oscf}(t,\vec{\ell})$ the scheme of \S\,\ref{sect:7} can be
used for evaluation of $\MF{J}(\vec{\ell},\vec{x})$ without any amendments or
restrictions.

Consider now application of the scheme developed in \S\,\ref{sect:3}, \ref{sect:5} to
computation of \eqref{eq:Jldef}. Similarly to \eqref{eq:Idef} we write
\[
 \MF{J}(\vec{\ell},\vec{x})=
 \int_0^1 \MF{\oscf_*}(\tau,\vec{\ell}) \frac{\E{\MF{\oscf_*}(\tau,\vec{x})}}{(1-\tau)^2}\,\D{}t,
\]
where again $\MF{\oscf_*}(\tau,\vec{x})=\MF{\oscf}(\tau/(1-\tau),\vec{x})$. Then
we use \eqref{eq:Levin0} and find that
\begin{equation}
 \MF{J}(\vec{\ell},\vec{x})=\frac{\MF{\Psi_{\vec{\ell},\vec{x}}}(\tau)}
 {(1-\tau)^2}\E{\MF{\oscf_*}(\tau,\vec{x})}\Bigm|_{\tau=0}^{\tau=1},
\label{eq:JviaPsi}
\end{equation}
if for the given $\bm{\ell}$ and $\vec{x}$ the function
$\MF{\Psi_{\vec{\ell},\vec{x}}}(\tau)$ satisfies the equation
\[
  \MF{\Psi'_{\vec{\ell},\vec{x}}}(\tau)+
  \frac{2\MF{\Psi_{\vec{\ell},\vec{x}}}(\tau)}{1-\tau}+
  \MF{\pd{\tau}\oscf_*}(\tau,\vec{x})\MF{\Psi_{\vec{\ell},\vec{x}}}(\tau)=
  \MF{\oscf_*}(\tau,\vec{\ell}).
\]

Since $\MF{\oscf_*}(\tau,\vec{\ell})=\MF{O}((1-\tau)^{-2})$ as $\tau\to1$, it is
convenient to define
$\MF{\Psi_{\vec{\ell},\vec{x}}}(\tau)=(1-\tau)\MF{\Phi_{\vec{\ell},\vec{x}}}(\tau)$.
Then, $\MF{\Phi_{\vec{\ell},\vec{x}}}(\tau)$ satisfies the following
differential equation
\begin{equation}
  (1-\tau)^3\MF{\Phi'_{\vec{\ell},\vec{x}}}(\tau)+\bigl[(1-\tau)^2+\MF{\sigma}(\tau,\vec{x})\bigr]\MF{\Phi_{\vec{\ell},\vec{x}}}(\tau)=
  \MF{\oscf_\circ}(\tau,\vec{\ell}),
\label{eq:der*}
\end{equation}
where $\sigma$ is defined by \eqref{eq:sigmadef} and
$\MF{\oscf_\circ}(\tau,\vec{\ell})=(1-\tau)^2\MF{\oscf_*}(\tau,\vec{\ell})$.

Analogously to \eqref{eq:Phi_gen}, \eqref{eq:Phi0def}, assuming $y<0$ we can find the
bounded on $[0,1]$ solution of \eqref{eq:der*}
\begin{equation*}
 \MF{\Phi_{\vec{\ell},\vec{x}}}(\tau)=(\tau-1)\E{\MF{\Lambda}(\tau,\vec{x})}\int_\tau^1
\frac{\MF{\oscf_\circ}(\tau,\vec{\ell})\E{-\Lambda(\theta,\vec{x})}}
{(\theta-1)^4}\,\D\theta,
\end{equation*}
where $\MF{\Lambda}(\tau,\vec{x})$ is defined by \eqref{eq:Lambdadef}. We write
\begin{equation*}
\MF{\Phi_{\vec{\ell},\vec{x}}}(\tau)=(\tau-1)\E{\MF{\Lambda}(\tau,\vec{x})}
 \int_\tau^1\frac{\E{-\frac{\gamma_2}{(\theta-1)^2}-\frac{\gamma_1}{\theta-1}-
 \gamma_0}}{(\theta-1)^4}\left[\MF{c^\ast_{0}}(\vec{\ell}) +
 \MF{\varrho}(\theta,\vec{\ell},\vec{x})\right]\,\D\theta,
\end{equation*}
where
$\MF{c^\ast_{0}}(\vec{\ell})=\MF{\oscf_\circ}(1,\vec{\ell})=\ell_2+\ii\ell_3$
and for $\theta\in[0,1]$, $\MF{|\varrho}(\theta,\vec{\ell},\vec{x})|\leq
\MF{C^\ast}(\vec{\ell},\vec{x})|\theta-1|$, $\MF{C^\ast}(\vec{\ell},\vec{x})$
is constant in $\theta$. Thus, we find
\begin{gather}
\MF{\Phi_{\vec{\ell},\vec{x}}}(\tau)= (\tau-1)\E{\MF{\Lambda}(\tau,\vec{x})}
\left\{ \MF{c^\ast_{0}}(\vec{\ell}) \MF{L^\ast_{0}}(\tau,\vec{x}) +
\MF{\check{L}}(\tau,\vec{\ell},\vec{x})\right\},\quad\mbox{where}
\label{eq:der*2.25--}
\\
\MF{L^\ast_{0}}(\tau,\vec{x})=\int_\tau^1(\theta-1)^{-4}\E{-\frac{\gamma_2}{(\theta-1)^2}
 -\frac{\gamma_1}{\theta-1}-\gamma_0}\,\D\theta,
\quad
\MF{\check{L}}(\tau,\vec{\ell},\vec{x})=\int_\tau^1
\frac{\MF{\varrho}(\theta,\vec{\ell},\vec{x})}{(\theta-1)^{4}}
\E{-\frac{\gamma_2}{(\theta-1)^2}
 -\frac{\gamma_1}{\theta-1}-\gamma_0}
\,\D\theta.\notag
\end{gather}
Further we obtain (cf.\ \eqref{eq:*2.5--})
\begin{equation}
  \MF{L^\ast_{0}}(\tau,\vec{x}) = \frac{\sqrt{\pi}}{8\gamma_2^{5/2}}(\gamma_1^2+2\gamma_2)\E{-\gamma_0+\frac{\gamma_1^2}{4\gamma_2}}
 \MF{\mathrm{erfc}}\left(\frac{\sqrt{\gamma_2}}{1-\tau}-\frac{\gamma_1}{2\sqrt{\gamma_2}}\right)+
\frac{1}{4\gamma_2^2}\E{-\frac{\gamma_2}{(\tau-1)^2}-\frac{\gamma_1}{\tau-1}-\gamma_0}
\left(\gamma_1+\frac{2\gamma_2}{1-\tau}\right).
\label{eq:der*2.5--}
\end{equation}

Define the function
\begin{equation}
 \lambda(\tau,y)=\frac{\sqrt{\pi}}{4|y|^{3/2}}(1+2|y|)\E{|y|}
 \MF{\mathrm{erfc}}\left(\frac{\sqrt{|y|}\tau}{1-\tau}\right)+
\frac{1}{2|y|}\E{\frac{y}{(\tau-1)^2}+\frac{2y}{\tau-1}}
\left(1+\frac{1}{1-\tau}\right)
\label{eq:derldef}
\end{equation}
(the expression is obtained by replacing $\gamma_0$, $\gamma_1$, $\gamma_2$ by
by their real parts in the right-hand side of \eqref{eq:der*2.5--}). Then, for
$\tau\leq1$ we find $\MF{|\check{L}}(\tau,\vec{\ell},\vec{x})|\leq
\MF{C^\ast}(\vec{\ell},\vec{x})(1-\tau)\lambda(\tau,y)$ and, thus, in view of
\eqref{eq:derldef} and \cite[7.1.23]{AS},  we have as $\tau\to1-0$:
\begin{equation}
\left|(\tau-1)\E{\MF{\Lambda}(\tau,\vec{x})}\MF{\check{L}}(\tau,\vec{\ell},\vec{x})\right|=\MF{O}(\tau-1).
\label{eq:derests}
\end{equation}
Using \eqref{eq:der*2.25--}, \eqref{eq:der*2.5--}, \eqref{eq:derests}, and
\cite[7.1.23]{AS} we obtain the representation
\begin{gather}
\MF{\Phi_{\vec{\ell},\vec{x}}}(\tau)=\MF{\Phi_{\vec{\ell},\vec{x}}}(1)
 +\MF{O}(\tau-1),\quad\mbox{where}
\notag\\
 \MF{\Phi_{\vec{\ell},\vec{x}}}(1):=\lim_{\tau\to1-0}\MF{\Phi_{\vec{\ell},\vec{x}}}(\tau)=
\frac{\ell_2+\ii\ell_3}{2(y+\ii z)}.
\label{eq:derPhi0}
\end{gather}

Finally, for $y<0$ by \eqref{eq:JviaPsi} we have
$\MF{J}(\vec{\ell},\vec{x})=-\MF{\Phi_{\vec{\ell},\vec{x}}}(0)\E{y+\ii x}$. The
numerical solution to \eqref{eq:der*}, \eqref{eq:derPhi0} can be sought using
the scheme described in \S\,\ref{sect:5}. When we seek the approximation of
$\MF{\Phi_{\vec{\ell},\vec{x}}}(\tau)$ in the form \eqref{eq:PhiM_Tm_Mod} we can
define the first term of the expansion as follows:
\[
 \MF{\skew{5}\Hat{\phi}_{\vec{\ell},\vec{x}}}(\tau)=\MF{c^\ast_{0}}(\vec{\ell})
 (\tau-1)\E{\frac{\gamma_2}{(\tau-1)^2}+\frac{\gamma_1}{\tau-1}+\gamma_0}
 \MF{L^\ast_{0}}(\tau,\vec{x})
\]
(see \eqref{eq:der*2.25--}). It is also straightforward to compute
$\MF{\mathcal{L}\skew{5}\Hat{\phi}_{\vec{\ell},\vec{x}}}(\tau)$ analogously to
\eqref{eq:l0phidef}.

\section{Conclusion}

In this paper we deal with the problem of evaluation of the Green's function for
the classical linear ship-wave problem describing forward motion of bodies in
unbounded heavy fluid having a free surface. Of interest for us is the wavelike
(often referred to as `single integral') term $\MF{I_\infty}(x,y,z)$ which
represents the dominating in the far-field, oscillatory part of Green's
function. The integral $\MF{I_\infty}(x,y,z)$ is expressed \eqref{eq:defI_inf}
in terms of the integral $\MF{I}(x,y,z)$. Our purpose is to elaborate accurate
and fast computation techniques to approximate $\MF{I}(x,y,z)$ and its
derivatives. At that, the main difficulty is due to the presence in the
integrand of two oscillating factors of different nature and the infinite
interval of integration. The oscillating integrand can have stationary points
and there is a difficult limiting case~---~the track of the source moving in the
free surface is a line of essential singularities of the considered integral.

First, by using the ideas of Levin \cite{Levin82,Levin97} we reduce evaluation
of the integral to solution of an ordinary differential equation on the interval
$[0,1]$. We prove that the equation has one bounded solution,
whose value at the right end of the interval is known, while the value at the
left end, up to a known factor, coincides with the sought value of
$\MF{I}(x,y,z)$. To find the solution of the differential equation numerically
we develop an algorithm based on usage of Lagrange interpolating polynomial in a
barycentric form and collocation of the equation on a set of Chebyshev points.
Another representation of the solution to the differential equation consists of
the polynomial and a term arising from an asymptotic analysis. An estimate for
the residual of solution is suggested and numerically tested to demonstrate its
reliability. It is notable that, unlike the standard Levin's collocation based
on expansion in Chebyshev polynomials, the suggested scheme demonstrates
numerical stability.

Secondly, we develop an alternative numerical algorithm based on the
Clenshaw--Curtis quadrature \cite{ClCu} and involving transformation of the
integral path using the steepest descent method. The advantages of the
quadrature rule are its fast convergence, simple and effective computation of
weights (even for very large number of nodes), excellent numerical stability.
Relying upon the properties and taking into account the simplicity of the
integrand, we suggest to apply the quadrature in `brute force' manner,
increasing the number of its nodes (doubled at each step) until some last values
of the sequence of approximations become closer to each other with the given
tolerance.

These two alternative methods are tested numerically and compared for a wide
variety of parameters, with special attention to the accuracy and efficacy. The
experiments show that both algorithms are reliable, compatible in speed and
much faster than standard solvers being applied to the Levin's differential
equation. The algorithm based on Levin's equation and barycentric Lagrange
interpolation is somewhat faster when the order of interpolating polynomial,
needed for achieving the given accuracy, is small. At the same time, the
algorithm based on the Clenshaw--Curtis quadrature is more convenient to
evaluate $\MF{I}(x,y,z)$ to the given accuracy and the algorithm works better in
the most difficult for the numerical integration zone near the track of the
source moving close to the free surface.

Finally, in the present work we discuss application of the suggested methods to
computation of derivatives of $\MF{I}(x,y,z)$ and show that their generalization
is rather straightforward.

\end{document}